\documentclass[12pt]{article}
\usepackage{color}
\usepackage{amscd,amsfonts,amssymb,amscd,amsmath,latexsym}

\usepackage[hypertex,backref]{hyperref}
\usepackage[all]{xy}
\usepackage{mathrsfs}

\title{Uniqueness of smooth extensions of generalized cohomology theories}
\author{Ulrich Bunke\thanks{NWF I - Mathematik,
Universit{\"a}t Regensburg,
93040 Regensburg,
GERMANY, ulrich.bunke@mathematik.uni-regensburg.de}  
 and Thomas Schick\thanks{Mathematisches Institut, Georg-August-Universit\"at
   G\"ottingen, 37073 G\"ottingen, Bunsenstr. 3-5,
   schick@uni-math.gwdg.de\protect\\ Thomas Schick was partly partially funded by the Courant Research Center "`Higher order structures in Mathematics"'
within the German initiative of excellence} }

\newtheorem{theorem}{Theorem}[section] 
\newtheorem{prop}[theorem]{Proposition}
\newtheorem{lem}[theorem]{Lemma}
\newtheorem{ddd}[theorem]{Definition}
\newtheorem{kor}[theorem]{Corollary}
\newtheorem{ass}[theorem]{Assumption}

\newcommand{\Rham}{{\tt Rham}}

\newcommand{\bMU}{{\mathbf{MU}}}
\newcommand{\tT}{{\tt T}}
\newcommand{\tB}{{\tt B}}

\newcommand{\hocolim}{{\tt hocolim\:}}

\newcommand{\colim}{{\tt colim}}

\newcommand{\Z}{\mathbb{Z}}

\newcommand{\bh}{{\mathbf{h}}}
\newcommand{\bE}{{\bf E}}
\newcommand{\bM}{{\bf M}}

\newcommand{\bU}{{\mathbf{U}}}

\newcommand{\proof}{{\it Proof.$\:\:\:\:$}}
\newcommand{\ori}{{\mathrm{ori}}}

\newcommand{\R}{\mathbb{R}}

\newcommand{\Q}{\mathbb{Q}}

\newcommand{\bF}{{\bf F}}

\newcommand{\tZ}{{\tt Z}}
\newcommand{\tV}{{\tt V}}

\newcommand{\bK}{{\bf K}}
\newcommand{\bo}{{\bf o}}

\newcommand{\bS}{\mathbf{S}}

\newcommand{\cE}{\mathcal{E}}

\newcommand{\cK}{\mathcal{K}}

\newcommand{\Hom}{{\tt Hom}}

\newcommand{\im}{{\tt im}}

\newcommand{\inter}{{\tt int}}

\newcommand{\coker}{{\tt coker}}
\newcommand{\id}{{\tt id}}

\newcommand{\nat}{\mathbb{N}}

\def\imath{{i}}

\def\hB{\hspace*{\fill}$\Box$ \newline\noindent}

\def\hB{\hspace*{\fill}$\Box$ \\[0.5cm]\noindent}

\newcommand{\pr}{{\tt pr}}

\newcommand{\ch}{{\mathbf{ch}}}

\newcommand{\Ab}{{\mathrm{Ab}}}

\begin{document}\maketitle
\begin{abstract}
We provide an axiomatic framework for the study of smooth extensions of
generalized cohomology theories. Our main results are about  the uniqeness of
smooth extensions,  and the identification of the flat theory with the
 associated cohomology theory with $\R/\Z$-coefficients. 

\textcolor{black}{In particular, we show that there is a unique smooth extension of K-theory and
of MU-cobordism with a unique multiplication, and that the flat theory in
these cases is naturally isomorphic to the homotopy theorist's version of the
cohomology theory with $\R/\Z$-coefficients. For this we only require a small
set of natural compatibility conditions.}
\end{abstract}

\tableofcontents

\section{Axioms}

A smooth extension of a generalized cohomology theory $E$ is a refinement $\hat E$ of the restriction of  $E$ to the category of smooth  manifolds. The functor $\hat E$ is no longer homotopy invariant. A class $\hat x\in \hat E(M)$ which refines the underlying topological class $I(\hat x)\in E^*(M)$ contains the information about a closed differential form
 $R(\hat x)\in \Omega_{cl}^*(M,E^*\otimes_\Z\R)$ which represents the image of $I(\hat x)$ under the natural map $\ch\colon E^*(M)\to H^*(M;E^*\otimes_\Z\R)$. The deviation of $\hat E$ from homotopy invariance is described by a homotopy formula (Lemma \ref{udqwdqwdqw1}). 
Let $\hat x\in \hat E^{*+1}([0,1]\times M)$ and $f_0,f_1\colon M\to [0,1]\times M$ be the inclusions of the endpoints. Then 
\begin{equation}\label{udqwdqwdqw1}
f_1^*\hat x-f_0^*\hat x=a(\int_{[0,1]\times M/M} R(\hat x))\ ,
\end{equation}
where $a$ is the natural transformation \ref{ddd1}.A.4.

A typical and motivating example is the smooth extension $\widehat{H\Z}^*$ of
integral cohomology. The group $\widehat{H\Z}^2(M)$ can be identified with the
group of isomorphism classes $[L,h^L,\nabla^L]$ of  hermitean line bundles  on
$M$ with unitary connections with the tensor product as the group operation.
We have
$I([L,h^L,\nabla^L])=c_1(L)\in H\Z^2(M)$, the first Chern class of $L$, and
$R([L,h^L,\nabla^L])=-\frac{1}{2\pi i}R^{\nabla^L}\in \Omega^2_{cl}(M)$, the first Chern form.
Unlike the first Chern class  $c_1(L)\in H\Z^2(M)$, the class
$[L,h^L,\nabla^L]\in \widehat{H\Z}^2(M)$ captures secondary information, e.g.\   the holonomy of $\nabla^L$ which might be non-trivial even if $L$ is trivial and $\nabla^L$ is flat. Refined characteristic classes for flat bundles
were one of the  first motivations for the introduction of
$\widehat{H\Z}$ in \cite{MR827262}.

The space of hermitean line bundles  with unitary connections is the configuration space of Maxwell field theory, i.e. the gauge theory with structure group $U(1)$. 
In this field theory the field strength is a closed two-form
which satisfies the   following quantization condition: The 
 integral of the field strength over cycles is required to be integral. In the past decade the  discussion of models of string theory with branes lead to field theories with $p$-form field strength. Furthermore, the quantization conditions motivated the consideration of  underlying cohomology theories
different from ordinary cohomology theory like $K$-theory, see e.g.  \cite{freed-2000}, \cite{freed-2000-0005}, \cite{moore-2000-0005}. 

The use of smoothly extended cohomology groups as configuration spaces in field theories, the topological considerations in \cite{MR2192936}, and further developments on secondary invariants (see e.g. \cite{bunke-2002} and the literature cited therein) lead to the development of this circle of ideas to a mathematical theory. The present paper contributes to this theory by presenting axioms for smooth extensions
and showing that they imply uniqeness results in many interesting cases.

We consider a generalized cohomology theory $E$ represented by a spectrum $\bE$.
It gives rise to the $\Z$-graded abelian group $E^*:=E^*(*)=\pi_{-*}\bE$, 
and we define the $\Z$-graded $\R$-vector space $\tV:=E^*\otimes_\Z \R$. For a smooth manifold $M$ we define
$\Omega^*(M,\tV):=C^\infty(M,\Lambda^*T^*M\otimes_\R \tV)$ with the $\Z$-grading by the total degree. To be more precise, in the case of infinite-dimensional $\tV^n$ we topologize $\tV^n$ as a colimit of its finite-dimensional subspaces with their canonical real vector space topologies.
Locally an element of $\Omega^*(M,\tV^n)$ can then be written
as a finite sum $\sum_{j} \omega_j\otimes v_j$ for collections of forms $\omega_j\in \Omega^*(M)$ and elements $v_j\in \tV^n$.
We let $d\colon \Omega^*(M,\tV)\to \Omega^{*+1}(M,\tV)$ be the de Rham differential, and we write
$\Omega_{cl}^{*}(M,\tV):=\ker(d\colon \Omega^*(M,\tV)\to \Omega^{*+1}(M,\tV))$ for the subspace of closed forms.  We identify
$H^*(M;\tV)$ with the singular cohomology of $M$ with coefficients in $\tV$. Integration over simplices induces the natural transformation
$$\Rham\colon \Omega_{cl}^{*}(M,\tV)\to H^*(M;\tV)\ .$$
It
induces an isomorphism $H^*_{dR}(M;\tV)\stackrel{\sim}{\to}
H^*(M;\tV)$. Furthermore, there
is a canonical natural transformation of cohomology theories
$$\ch\colon E^*(X)\to H^*(X;\tV)\ .$$

\begin{ddd}\label{ddd1}
A \textbf{smooth extension} of the generalized cohomology theory $E$ is a quadruple
$(\hat E,R,I,a)$, where
\begin{enumerate}
\item[A.1] $\hat E$ is a contravariant functor from the category of smooth manifolds to $\Z$-graded abelian groups. Sometimes we will consider a version defined only on the category of compact manifolds (possibly with boundary). 
\item[A.2] $R$ is a natural transformation of  $\Z$-graded abelian group-valued functors
$$R\colon \hat E^*(M)\to \Omega^*_{cl}(M,\tV)\ .$$
\item[A.3] $I$ is a natural transformation of  $\Z$-graded abelian group-valued functors
$$I\colon \hat E^*(M)\to E^*(M)\ .$$
\item[A.4] $a$ is a natural transformation of  $\Z$-graded abelian group-valued functors
$$a\colon \Omega^{*-1}(M,\tV)/\im(d)\to \hat E^{*}(M)\ .$$
\end{enumerate}
These objects have to satisfy the following relations:
\begin{enumerate}
\item[R.1]  
$R\circ a=d$
 \item[R.2]  For all manifolds $M$ the diagram $$\xymatrix{\hat E^*(M)\ar[d]^{I}\ar[r]^{R}&\Omega_{cl}^*(M,\tV)\ar[d]^{\Rham}\\E^*(M)\ar[r]^{\ch}&H^*(M;\tV)}$$
commutes.
\item[R.3] For all manifolds $M$ the sequence \begin{equation}\label{ddd13}E^{*-1}(M)\stackrel{\ch}{\to} \Omega^{*-1}(M)/\im(d)\stackrel{a}{\to} \hat E^*(M)\stackrel{I}{\to} E^*(M)\to 0
\end{equation}
is exact.
\end{enumerate}
\end{ddd}

We now consider two smooth extensions
$(\hat E,R,I,a)$ and $(\hat E^\prime,R^\prime,I^\prime,a^\prime)$ of $E$. 
\begin{ddd}
A \textbf{natural transformation of smooth extensions} is a natural transformation of $\Z$-graded abelian group valued functors
$\Phi\colon \hat E^*\to \hat E^{\prime *}$ such that the following diagram
commutes for every manifold $M$:
$$\xymatrix{\Omega^{*-1}(M,\tV)\ar[r]^a\ar@{=}[d]&\hat E^*(M)\ar[d]^{\Phi}\ar[r]^{I}\ar@/^1cm/[rr]^R&E^*(M)\ar@{=}[d]&\Omega^*_{cl}(M,E)\ar@{=}[d]\\\Omega^{*-1}(M,\tV)\ar[r]^{a^\prime}& \hat E^{\prime*}(M)\ar[r]^{I^\prime}\ar@/_1cm/[rr]^{R^\prime}&E^*(M)&\Omega^*_{cl}(M,E)}
\ .$$
\end{ddd}

We consider the inclusion of the base point $*\to S^1$. It induces an embedding
$i\colon M\to S^1\times M$ for every manifold $M$. Let $p\colon S^1\times M\to M$ be the projection onto the second factor. Since $p\circ i=\id_{M}$ we get splittings
$$\hat E^*(S^1\times M)\cong \im(p^*)\oplus \ker(i^*) ,\quad   E^*(S^1\times M)\cong \im(p^*)\oplus \ker(i^*)\ .$$ Let $\Sigma M_+$ be the suspension which is a space, not a manifold. There is a natural
projection $q\colon S^1\times M\to \Sigma M_+$ of spaces. It induces an isomorphism
$$q^*\colon \tilde E^*(\Sigma M_+)\stackrel{\sim}{\to} \ker(i^*)\subseteq E^*(S^1\times M)\ .$$
Furthermore, there is the suspension isomorphism 
$$\sigma\colon E^{*-1}(M)\stackrel{\sim}{\to} \tilde E^*(\Sigma M_+)\ .$$
Composing these isomorphisms with the projection onto $\ker(i)$ we get the integration map
$$\int\colon E^{*+1}(S^1\times M)\stackrel{\pr}{\to} \ker(i^*)\stackrel{(q^{*})^{-1}}{\to}  \tilde E^{*+1}(\Sigma M_+)\stackrel{\sigma^{-1}}{\to} E^{*}(M)$$
for the generalized cohomology theory $E$.

We introduce the notation $SF(M):=F(S^1\times M)$ for a functor $F$ defined on manifolds.
The integration $$\int\colon SE^{*+1}\to E^*$$ just defined is complemented by an integration map
$$\int\colon S\Omega^{*+1}(\dots,\tV)\to \Omega^{*}(\dots,\tV)$$ for differential forms which preserves the image and kernel of $d$.

\begin{ddd}\label{ddd4}
A smooth extension \textbf{with integration} of $E$ is a quintuple $(\hat E,R,I,a,\int)$, where $(\hat E,R,I,a)$ is a smooth extension of $E$ and
$\int$ is a natural transformation $$\int\colon S\hat E^{*+1}\to \hat E^*$$ such that
\begin{enumerate}
\item $\int \circ (t^*\times \id)^*=-\int$, where $t\colon S^1\to S^1$ is given by $t(z):=\bar z$.
\item \label{ddd41} $\int\circ p^*=0$ and
\item
the  diagram
 $$\xymatrix{S\Omega^{*}(M,\tV)\ar[r]^a\ar[d]^\int&S\hat E^{*+1}(M)\ar[d]^{\int}\ar[r]^{I}\ar@/^1cm/[rr]^R&SE^{*+1}(M)\ar[d]^\int&S\Omega^{*+1}_{cl}(M,E)\ar[d]^\int\\\Omega^{*-1}(M,\tV)\ar[r]^{a}& \hat E^{*}(M)\ar[r]^{I}\ar@/_1cm/[rr]^{R}&E^*(M)&\Omega^*_{cl}(M,E)}
$$ commutes for all manifolds $M$.\end{enumerate} 
\end{ddd}

We now consider two smooth extensions with integration
$(\hat E,R,I,a,\int)$ and $(\hat E^\prime,R^\prime,I^\prime,a^\prime,\int^\prime)$ of $E$. 
\begin{ddd}
A \textbf{natural transfomation  between two extensions with integration} of $E$
is a natural transfomation  between smooth extensions $\Phi\colon \hat E^*\to \hat E^{\prime,*}$
such that
$$\xymatrix{S\hat E^{*+1}(M)\ar[d]^\int\ar[r]^\Phi&S\hat E^{\prime*+1}(M)\ar[d]^{\int^\prime}\\\hat  E^*(M)\ar[r]^{\Phi}&\hat E^{\prime*}(M)}$$
commutes for all manifolds $M$.
\end{ddd}

Assume now that $E$ is a multiplicative cohomology theory. Then the  functor $E^*$ has values in graded commutative rings. In particular, $E^*$ is a $\Z$-graded ring, and $H^*(M;\tV )$ and $\Omega^*(M,\tV )$ are  $\Z$-graded  rings as well. 
In this case we can make the following definition.
\begin{ddd}
A \textbf{multiplicative smooth extension} is a smooth extension $(\hat E,R,I,a)$
such that $\hat E^*$ takes values in $\Z$-graded commutative rings, the transformations $R$ and $I$ are multiplicative, and the identity
$$x\cup a(\alpha)=a(R(x)\wedge \alpha)\:\:\footnote{Observe, that $R(x)\wedge d\alpha=d(R(x)\wedge \alpha)$ so that the right-hand side is well-defined.}\ ,\quad x\in \hat E^*(M)\ , \quad \alpha\in \Omega^{*}(M,\tV )/\im(d)$$
holds true.
\end{ddd}

For every generalized cohomology theory $E^*$ represented by a spectrum $\bE$
a smooth extension $(\hat E,R,I,a)$ exists by the construction of Hopkins-Singer \cite{MR2192936}.

The historically first example of a smooth extension was constructed by Cheeger and Simons in  \cite{MR827262} for ordinary integral cohomology $H\Z$ (and more generally for $H R$ for discrete subrings $R$ of $\R$). These extensions of ordinary cohomology  are multiplicative.
The classes in $\widehat{H R}(M)$ were realized in \cite{MR827262} as differential characters. 
By now there are various different constructions of the smooth extension of ordinary cohomology, e.g by sheaf theory \cite{MR1197353} (under the name smooth Deligne cohomology), using geometric cycles Gajer \cite{MR1423029}, cubical chains \cite{MR2179587}, or stratifold bordisms \cite{bks}.
With differential characters the integration over $S^1$ as in Definition \ref{ddd4}, but also for general proper submersions $p\colon M\to B$ is simple, but the product is complicated. 
In cochain  models both structures are involved, while in the sheaf-theoretic Deligne cohomology
the product is easy, but integration is complicated. In the stratifold bordism model both structures
are straight forward and explicit, and therefore this model is predestinated for the verification of the projection formula
$\int_p (y\cup p^*x)=(\int_p y)\cup x$, where $y\in \widehat{H R}^*(M)$, $x\in \widehat{H R}^*(B)$.

In view of the variety of  constructions of a smooth extension of ordinary cohomology it is a natural question whether all give equivalent results.
This has been answered by \cite{MR2365651}, though a set of slightly different axioms is used\footnote{In  \cite{MR2365651} the additional requirement is that the flat theory (Definition \ref{uiqwdqwdwqd}) is topological (Definition \ref{udiqdwqdwqd54545}).}.
Uniqueness also follows from the axioms stated above  by Theorem \ref{main1}
if one takes the integration into account. This has first been observed by
 Moritz Wiethaup (2006/2007). In both cases the smooth
extension is unique up to unique isomorphism. Moreover, we have uniqueness of
the product by   \cite[Thm. 1.3]{MR2365651} or Theorem \ref{main23}.

In \cite{bunke-2007} we give a geometric construction of smooth extensions of
bordism theories.
We developed the details in the case of complex bordism $MU$. The method also applies to other bordism theories, e.g. oriented bordism, $MSpin$ or $MSpin^c$-bordism or framed bordism $S$. In all these cases we obtain a multiplicative extension and a theory of integration for suitably oriented proper submersions. The particular importance of the case of complex bordism theory comes from the Landweber exact functor theorem \cite{MR0423332}. It allows to construct a multiplicative smooth extension for every complex oriented Landweber exact cohomology theory. Examples are complex $K$-theory and certain elliptic cohomology theories. 

Using methods of local index theory in \cite{bunke-20071} we have constructed a Dirac operator model of smooth $K$-theory which is again multiplicative and has a nice integration theory for smoothly $K$-oriented proper submersions. 

The presence of different constructions (at least two in the case of bordism theories \cite{MR2192936} and \cite{bunke-2007}, and three in the case of $K$-theory \cite{MR2192936},  \cite{bunke-20071}, \cite{bunke-2007}) raises again the question whether they are equivalent. Moreover, for applications to topology, e.g. constructions of secondary invariants, of particular importance is the identification of the associated flat theory with the corresponding $\R/\Z$-theory.  

To answer these questions is the main motivation and result of the present paper. Note that all these examples are rationally even (Definition \ref{uzdiqwdqwdwqd666}).
The examples constructed from $MU$ by the Landweber exact functor theorem are only defined on the category of compact manifolds. This is the reason for considering this case in the present paper, too.
Observe that the coefficients of a Landweber exact theory are torsion-free. Therefore a rationally even Landweber exact cohomology theory is even. This is exactly the additional assumption made in the statement of our uniqueness Theorems \ref{main1} and   \ref{main23} in order to cover smooth extensions which are only defined on compact manifolds.

Let us now formulate the main results of the present paper.
\begin{theorem}[Thm. \ref{main1}] 
Let $E$ be a rationally even generalized cohomology theory represented by a spectrum $\bE$.  Let
$(\hat E,R,I,a,\int)$ and $(\hat E^\prime,R^\prime,I^\prime,a^\prime,\int^\prime)$ be  two smooth extensions with integration. We assume that either the smooth extensions are defined on the category of all smooth manifolds and the coefficients $E^m$ are countably generated for all $m\in \Z$, or that  $E^m=0$ for all odd $m\in \Z$ and $E^m$ is finitely generated for all even $m\in \Z$.
Then there is a unique natural isomorphism
$$\Phi\colon \hat E\to \hat E^\prime$$ of smooth extensions with integration.
\end{theorem}

A multiplicative smooth extension of a rationally even cohomology theory has a canonical integration by Corollary \ref{uifqfqfqf}. 

\begin{theorem}[Cor. \ref{main2}] 
Assume that  $(\hat E,R,I,a)$ and  $(\hat E^\prime,R^\prime,I^\prime,a^\prime)$ are two multiplicative extensions of a rationally even generalized cohomology theory with countably generated coefficients such that either both are defined on the category of all smooth manifolds or
$E^*$ is even and degree-wise finitely generated. Then there is a unique natural isomorphism between these  smooth extensions preserving the canonical integration. This transformation is multiplicative. 
\end{theorem}
At the moment we have no feeling how important the condition of $E$ beeing rationally even is.
This theorem applies e.g. to multiplicative smooth extensions of ordinary cohomology, all the bordism theories listed above and
complex or real $K$-theory (for complex $K$-theory also to  the version defined on compact manifolds).

The flat theory
$$\hat E^*_{flat}(M):=\ker\left(R\colon \hat E^*(M)\to \Omega_{cl}^*(M,\tV )\right)$$
is a homotopy invariant functor on smooth manifolds with values $\Z$-graded abelian groups. 
\begin{theorem}[Thm. \ref{zuddwqdqwdqw}]
Assume that $E$ is rationally even with countably generated coefficients. If $(\hat E,R,I,a,\int)$ is a smooth extension of $E$ with integration  which is defined on all smooth manifolds (or alternatively, on all compact manifolds and $E^*$ is even and degree-wise finitely generated),
then there exists an isomorphism
$$\Phi_{flat}\colon \hat E^*_{flat}\stackrel{\sim}{\to} E\R/\Z^{*-1}\ .$$
\end{theorem}
For a more precise statement and the notation see Section \ref{uciacascc}.
This theorem implies  that the additional axiom in \cite{MR2365651}
follows from our axioms together with the presence of integration.
Theorem \ref{zuddwqdqwdqw} in particular states that the flat theory $\hat E^*_{flat}$ is a generalized
cohomology theory. The essential additional datum turning a homotopy invariant functor into a generalized
cohomology theory is the boundary operator of a long exact Mayer-Vietoris sequence.
Theorem \ref{zuddwqdqwdqw} is proven by a comparison with the Hopkins-Singer theory \cite{MR2192936}.
In Section \ref{e89wfoewfwefq} we give an independent construction of the Mayer-Vietoris sequence.
\begin{theorem}[Thm. \ref{udidowqdiwqdwqd}]\label{uiewefwef}If $(\hat E,R,I,a,\int)$ is a smooth extension of $E$ with integration, then  $\hat E^*_{flat}$
has a natural long exact Mayer-Vietoris sequence. Its restriction to compact manifolds
is equivalent to the restriction  to compact manifolds of a generalized
cohomology \textcolor{black}{theory} represented by a spectrum.
\end{theorem}
Note that Theorem \ref{uiewefwef} does not require any additional assumptions, but it also does not state that
$\hat E^*_{flat}$ is equivalent to  $E\R/\Z^{*-1}$. This equivalence can be derived again under additional assumptions, now independently from  \cite{MR2192936}, but at the cost of restricting to compact manifolds.
\begin{theorem}[Thm. \ref{udqiwdqwdqwddwqdqwd1}]\label{udqiwdqwdqwddwqdqwd}
If $(\hat E,R,I,a,\int)$ is a smooth extension of $E$ with integration and $E^*$ is torsion free,  
then  we have a natural equivalence $\hat E^*_{flat}\cong E\R/\Z^{*-1}$ of functors restricted to the category of compact manifolds.
\end{theorem}

\textit{Acknowledgement: The basic questions for the present paper have been the topic of a PhD project of Moritz  Wiethaup (G\"ottingen). Some basic ideas and first results are due to him. In this paper we work out in detail and further develop the results of fruitful mathematical discussions in G\"ottingen arround 2006/2007.}

\section{Approximation of spaces by manifolds}\label{uziedwqeqw}
\newcommand{\Top}{{\mathcal{T}op}}

%


The main technical problem of the present paper is the construction of a natural transformation between two
smooth extensions $\hat E^k,\hat E^{\prime k}$ of a generalized cohomology
$E^k$. This requires to define a homomorphism $\hat E^k(M)\to\hat E^{\prime
  k}(M)$  for every smooth manifold in a natural way. If the underlying
topological functor $E^k$  for fixed $k\in\Z$  would be represented by a manifold $\bE$, naturality of the construction could be obtained by making one universal choice $\hat E^k(\bE)\to\hat E^{\prime k}(\bE)$, only. Unfortunately, a generalized cohomology functor $E^k$
is rarely represented by a finite-dimensional manifold $\bE$. The idea for nevertheless cutting down the amount of choices in order to secure naturality is to approximate the classifying space $\bE$ by a sequence of manifolds. Since in some examples our smooth extensions are only defined on compact manifolds we take care of the case where such an approximation can be realized by compact manifolds.

Recall that a map $f\colon X\to Y$ between pointed topological spaces is called $n$-connected if $f_*\colon \pi_k(X)\to \pi_k(Y)$ is an isomorphism for $k< n$ and surjective for $k=n$.

\begin{prop}\label{duwidwdwd}
 Let $\bE$ be a connected pointed topological space.
If  $\bE$  is simply connected and $\pi_k(\bE)$ is finitely generated for all $k\ge 2$,
then there exist a sequence of compact pointed manifolds with boundary
 $(\cE_i)_{i\in \nat}$ together with pointed maps $\kappa_i\colon \cE_i\to \cE_{i+1}$, $x_i\colon \cE_{i}\to \bE$ for all  $i\in \nat$ such that
\begin{enumerate}
\item\label{pr00} $\cE_i$ is homotopy equivalent to an $i$-dimensional $CW$-complex,
\item\label{pr01} the map $x_i$ is   $i$-connected,
\item\label{pr02} $\kappa_i\colon \cE_i\hookrightarrow   \cE_{i+1}$ is an embedding of a submanifold,
\item\label{pr03} the diagram
$$\xymatrix{\cE_{i}\ar[dr]_{x_i}\ar[rr]^{\kappa_i}&&\cE_{i+1}\ar[dl]^{x_{i+1}}\\&\bE&}$$ commutes,
\item \label{la}for  all finite-dimensional  pointed $CW$-complexes $X$ the induced map
$$\colim ( [X,\cE_{i}])\to [X,\bE]$$ is an isomorphism.
\end{enumerate}
\end{prop}
\proof
We first construct a similar sequence $(E_i)_{i\in \nat }$ of connected finite \textcolor{black}{CW-}complexes together with maps
$y_i\colon E_i\to \bE$ and $\sigma_i\colon E_i\to E_{i+1}$ so that \begin{enumerate}
\item $E_i$ is an $i$-dimensional \textcolor{black}{CW-}complex,
\item  the map $y_i$ is  $i$-connected
\item\label{pr2} $\sigma_i\colon E_i\hookrightarrow   E_{i+1}$ is a cofibration,
\item\label{pr3} the diagram $$\xymatrix{E_{i}\ar[dr]_{y_i}\ar[rr]^{\sigma_i}&&E_{i+1}\ar[dl]^{y_{i+1}}\\&\bE&}$$ commutes. \end{enumerate}
We let $E_{1}$ be a point.
Assume that we have constructed $E_j$ for all $1\le j<i$ such that $\pi_{j}(E_{j})$ is  finitely generated.
Then we have a surjective map
$$y_{i-1,*}\colon \pi_{i-1}(E_{i-1})\to \pi_{i-1}(\bE)\ .$$

We claim that
$$\ker \left(y_{i-1,*}\colon \pi_{i-1}(E_{i-1})\to \pi_{i-1}(\bE)\right)$$ is finitely generated.
We have $\pi_{1}(E_{1})=0$.
For $j\ge 2$ we know that $\pi_1(E_j)=0$ by our inductive assumption. The
homotopy groups of a finite simply-connected $CW$-complex are finitely
generated (see e.g.~\cite{hatcherAT}).  Our kernel is finitely generated,
since it is a subgroup of a finitely generated abelian group. This finishes
the verification of the claim.

Let $(\gamma_\iota\colon S^{i-1}\to \bE_{i-1})_{\iota\in I}$ be a finite family of representatives of generators of the kernel of $y_{i-1,*}\colon \pi_{i-1}(E_{i-1})\to \pi_{i-1}(\bE)$.
  Then we define
$$\tilde E_i:=E_{i-1}\cup_{(\gamma_\iota)_{\iota\in I}}\bigsqcup_{\iota\in I}D^i\ .$$
The map $y_{i-1}$ has an extension $\tilde y_i\colon \tilde E_i\to \bE$.
Note that
$\tilde y_{i,*}\colon \pi_j(\tilde E_i)\to \pi_j(\bE)$ is now an isomorphism for $j<i$.


Let $(\theta_\lambda\colon S^i\to \bE)_{\lambda\in L}$ be a finite set of generators for
$\pi_i(\bE)$. We form
$$E_i:=\tilde E_i\vee \bigvee_{\lambda\in L} S^i\ ,$$ and we extend
$\tilde y_i$ to $y_i\colon E_i\to \bE$ using the maps $\theta_\lambda$.
Then
$y_{j,*}\colon \pi_j(E_i)\to \pi_j(\bE)$ is an isomorphism for $j< i$ and surjective
for $j=i$.

We let $\sigma_{i-1}\colon E_{i-1}\to E_i$ be the inclusion.

We now construct the desired sequence of manifolds together with homotopy equivalences
$z_i\colon \cE_i\to E_i$ such that 
$$\xymatrix{\cE_i\ar[d]^{z_i}\ar[r]^{\kappa_i}&\cE_{i+1}\ar[d]^{z_{i+1}}\\E_i\ar[r]^{\sigma_i}&E_{i+1}}$$
commutes.

We again start with a point $\cE_{1}:=*$.
Assume that we have constructed
$\cE_{i-1}$. We choose an embedding
$(\cE_{i-1},\partial \cE_{i-1})\hookrightarrow (\R^n_+,\R^{n-1})$,
where $\R_+:=\{(x_1,\dots,x_n)|x_n\ge 0\}$, $\R^{n-1}\subset \R^n$ is identified with the boundary $\{x_n=0\}$, and $n$ is sufficiently large. 
 We let $p\colon U\to \cE_{i-1}$ be the projection
from a tubular neighbourhood. We choose smooth maps
$\tilde \gamma_\iota\colon S^{i-1}\to \partial U$ so that
$z_{i-1}\circ p\circ\tilde \gamma_\iota\sim \gamma_\iota$.
If we take $n\ge 2i$, then we can assume after a homotopy that these extend to a smooth embedding \begin{equation}\label{eqq1}\sqcup_{\iota\in I} \tilde \gamma_\iota\colon \bigsqcup_{\iota\in I}  S^{i-1}\times D^{n-i}\to \partial U\setminus \R^{n-1}\ .\end{equation}  We use these maps  in order to define
$$\tilde \cE_i:=U\cup_{\sqcup_{\iota\in I} \tilde \gamma_\iota  } \bigsqcup_{\iota\in I}  D^i\times D^{n-i}$$ (smooth out corners).
The map $z_{i-1}$ has a natural  (up to homotopy)  extension to
$\tilde z_i\colon \tilde \cE_i\to \tilde E_i$ which is again a homotopy equivalence.
By a similar procedure we form the boundary connected sum 
$\cE_i:=\tilde \cE_i\bigsqcup \sqcup_{\lambda\in L} S^i\times D^{n-i}$ and extend $\tilde z_i$ to $z_i\colon \cE_i\to E_i$
which is again a homotopy equivalence.

The map $\kappa_{i-1}\colon \cE_{i-1}\to \cE_i$ is the inclusion.
Furthermore we let $x_i\colon \cE_{i}\to \bE$ be the  composition  $x_i\colon \cE_i\stackrel{z_i}{\to}E_i\stackrel{y_i}{\to} \bE$.  
By construction, our sequence satisfies \ref{pr00}..\ref{pr03}.

\begin{lem}
Property  \ref{la} holds true.
\end{lem}
We work in the model category structure on pointed topological spaces $\Top_*$ where
\begin{enumerate}
\item weak equivalences are  $\pi_*$-equivalences,
\item fibrations are Serre fibrations,
\item cofibrations are defined by the left lifting property.
\end{enumerate}
It is known that retracts of relative $CW$-extensions are cofibrations
\cite[Thm. 2.4.19]{MR1650134}.
We consider the poset $\nat$ as a category.
The diagram category $\Top_*^\nat$ has a model category structure with
\begin{enumerate}
\item level weak equivalences,
\item cofibrations  $(X_i)\to (Y_i)$ are characterized by the latching space condition, which in this case reduces to
the property that $X_i\sqcup_{X_{i-1}}Y_{i-1}\to Y_i$ is a cofibration for all $i\ge 1$,
\item fibrations are level fibrations.
\end{enumerate}
We refer to \cite[Ch. 5.1]{MR1650134} for details. It follows
that a system $(X_i)$ is cofibrant if all  maps $X_{i-1}\to X_i$ are cofibrations.

The homotopy colimit $\hocolim\colon \Top^\nat_*
\to \Top_*$ is the left derived functor of the colimit
$\colim\colon \Top_*^{\nat}\to \Top_*$.

By construction we have a weak equivalence
$\colim  (E_i)\to \bE$. 
Since the structure maps
$E_i\to E_{i+1}$ are all cofibrations (since they are $CW$-extensions) the whole system $(E_i)$ is cofibrant in $\Top_*^\nat$, and  there is a homotopy equivalence $\hocolim (E_i)\to \colim (E_i)$.
The map
$(\cE_i) \to (E_i)$  is a level equivalence.
If $X$ is a finite-dimensional $CW$-complex, then we have naturally induced isomorphisms
$$\colim ([X,\cE_i])\cong \colim([X,E_i])\stackrel{!}{\cong} [X,\hocolim (E_i)]\cong [X,\colim (E_i)]\cong [X,\bE]\ .$$
At the marked isomorphism we use the  cellular approximation theorem and that $X$ is finite-dimensional.
\hB

%
%
%
%
%
%
%
%

The property that $\bE$ is simply connected has been used  in order to conclude  from the fact that  $E_i$ is a finite simply-connected  $CW$-complex that
$\pi_i(E_i)$ is finitely generated. Finiteness is needed since we want to approximate by compact manifolds. If we allow non-compact manifolds, then essentially the same proof
gives the following:

\begin{prop}\label{duwidwdwd1}
 Let $\bE$  be a topological space with countably many connected components such that
the groups $\pi_k(\bE,x)$ are countably generated for all $k\ge 1$ and $x\in \bE$.
Then there exist a sequence of pointed manifolds 
 $(\cE_i)_{i\in \nat_0}$ together with pointed maps $\kappa_i\colon \cE_i\to \cE_{i+1}$, $x_i\colon \cE_{i}\to \bE$ for all  $i\in \nat_0$ such that
\begin{enumerate}
\item\label{pr000} $\cE_i$ is homotopy equivalent to an $i$-dimensional \textcolor{black}{CW-}complex,
\item\label{pr010} the map $x_i$ is  $(i-1)$-connected,
\item\label{pr020} $\kappa_i\colon \cE_i\hookrightarrow   \cE_{i+1}$ is an embedding of a submanifold,
\item\label{pr030} the diagram
$$\xymatrix{\cE_{i}\ar[dr]_{x_i}\ar[rr]^{\kappa_i}&&\cE_{i+1}\ar[dl]^{x_{i+1}}\\&\bE&}$$ commutes\ ,
\item \label{la0} for  all finite-dimensional pointed $CW$-complexes $X$ the induced map
$$\colim  ([X,\cE_{i}])\to [X,\bE]$$ is an isomorphism.
\end{enumerate}
\end{prop}
\proof
If $\bE$ is not connected, then we can approximate the countably many components of $\bE$
separately. 
In the connected case
we construct the $CW$-approximations $(E_i,\sigma_i,y_i)$ as before, but
starting the induction with $i=0$ and $E_0:=*$.
Then we proceed with the construction of the family $(\cE_i,\kappa_i,x_i)$. The only modification is as follows.
In order to find the embedding (\ref{eqq1}) we compose the embedding
$(\cE_{i-1},\partial \cE_{i-1})\to (\R^n_+,\R^{n-1})$ with an embedding
$(\R^n_+,\R^{n-1})\to (\R^{n+1}_+,\R^{n})$, and we use the extra dimension in order
to separate the images of the $\tilde \gamma_{\iota}$. If we want the
manifolds $\cE_i$ without boundary we can just \textcolor{black}{remove} the boundary
\textcolor{black}{without changing the homotopy type}.
\hB


In the following we discuss further properties of the approximations found in Proposition \ref{duwidwdwd}
or \ref{duwidwdwd1}. The first is a certain Mittag-Leffler condition.
For $j>i$ let $\kappa_i^{j}:=\kappa_{j-1}\circ \dots \circ \kappa_i\colon \cE_i\to \cE_j$. We set
$\kappa_i^i:=\id_{\cE_i}$. 
Let $\tV :=\bigoplus_{n\in \Z} \tV^n$ be some $\Z$-graded abelian group and let
$$H^*(X;\tV ):=\prod_{s+t=*}H^{s}(X;\tV^t)$$ denote the ordinary cohomology of the space $X$ with coefficients in $\tV $. 
\begin{lem}\label{cikdcdc}
For all $i\in \nat_0$ and all $l\ge  1$ 
we have the equality of subgroups of $H^k(\cE_{i};\tV )$
$$\kappa_i^{i+l,*}H^k(\cE_{i+l};\tV ) = \kappa_i^{*}H^k(\cE_{i+1};\tV )\ .$$
\end{lem}
\proof
Note that an $i$-connected map
$f\colon X\to Y$ induces an isomorphism in cohomology
$f^*\colon H^s(Y;G)\to H^s(X;G)$ for $s\le i-1$, and an injection
for $s=i$, where $G$ is an arbitrary abelian group. 
Since $\cE_i$ is $i$-dimensional we have
$$H^k(\cE_i;\tV )\cong \bigoplus_{s+t=k,s\le i}H^s(\cE_i;\tV^t)\ .$$
Note that $\kappa_i^{i+l}$ is $i$-connected for all $l\ge 0$.
We therefore have for $s\le i$, $j>i$, $l\ge 1$
$H^s(\cE_j;\tV^t)\cong \kappa_j^{j+l\textcolor{black}{,*}} H^s(\cE_{j+l};\tV^t)$.
This implies for all $l\ge 1$ that
$$\kappa_i^{*}H^k(\cE_{i+1};\tV )\cong
\kappa_i^* \circ \kappa_{i+1}^{i+l,*}H^k(\cE_{i+l};\tV )\cong
\kappa_i^{i+l,*} H^k(\cE_{i+l};\tV )\ .$$
\hB


\begin{prop}\label{fudif}Let $V$ be as above and
  $\bE$, $(\cE_i)_{i\ge 0}$ as in Proposition \ref{duwidwdwd1}.
 We consider a class $u\in H^*(\bE;\tV )$.
There exists a sequence of forms  $\omega_i\in \Omega_{cl}^*(\cE_i,\tV )$
such that $\Rham(\omega_i)=x_i^*u$ and $\kappa_i^*\omega_{i+1}=\omega_i$ for all $i\ge 0$.
\end{prop}
\proof
We construct $\omega_i$ inductively.
Assume that $\omega_j$ has been constructed for all $j\le i$. 
There exists  $\tilde \omega_{i+1}\in \Omega_{cl}^*(\cE_{i+1},\tV )$
such that $\Rham(\tilde \omega_{i+1})= x_{i+1}^*u$. Since $\kappa_i^*\Rham(\tilde \omega_{i+1})=\kappa_i^*x_{i+1}^*u=x_i^*u$ there 
exists furthermore a form $\alpha\in \Omega^{*-1}(\cE_{i},\tV )$
such that
$\kappa_i^*\tilde\omega_{i+1}+d\alpha=\omega_i$.
Since $\kappa_i$ is a closed embedding of a submanifold there exists an extension
$\tilde \alpha\in \Omega^{*-1}(\cE_{i+1},\tV )$ such that $\kappa_i^*\tilde \alpha= \alpha$.
We then define $\omega_{i+1}=\tilde \omega_{i+1}+d\tilde \alpha$.\hB

Let $(\hat E,R,I,a)$ be a smooth extension of a cohomology theory $E$. 
We  assume that this smooth extension is defined on all smooth manifold, or we assume
that our approximation $(\cE_i,x_i,\kappa_i)$ of $\bE$ is by compact manifolds. We consider the $\Z$-graded vector space $\tV :=E^*\otimes \R$. 
By Proposition \ref{fudif} we choose a sequence of closed forms $\omega_i\in \Omega^*_{cl}(\cE_i,\tV )$
such that $\Rham(\omega_i)=\ch(x_i^*u)$ and $\kappa_i^*\omega_{i+1}=\omega_i$ for all $i\ge 0$.
\begin{prop}\label{89fwefwefewfwe}
There exists a sequence $\hat u_i\in \hat E^*(\cE_i)$ such that $R(\hat u_i)=\omega_i$,
$I(\hat u_i)=x_i^*u$ and $\kappa_i^*\hat u_{i+1}=\hat u_i$ for all $i\ge 0$.
\end{prop}
\proof
%
%
First we  choose for each $i$ independently a class $ \tilde u_i\in \hat E^*(\cE_i)$ such that
$R(\tilde u_i)=\omega_i$ and $I(\tilde u_i)=x_i^*u$.
Note that $\tilde u_i$ is unique up to addition of a term $a(\alpha_i)$ for $\alpha_i\in H^{*-1}(\cE_i,\tV )$. 

We now argue by induction over $i$ using the Mittag-Leffler condition \ref{cikdcdc} that we can modify
our choice by $\hat u_i:=\tilde u_i+a(\alpha_i)$ so that $\kappa_i^*\hat u_{i+1}=\hat u_i$.

Let us assume by induction that we have made the choice of $\hat u_i$, and that we have already
chosen $\hat u_{i+1}^\prime\in \hat E^*(\cE_{i+1})$ such that $\kappa_i^*(\hat u_{i+1}^\prime)=\hat u_{i}$.
 Then $\kappa_{i+1}^*\tilde  u_{i+2}-\hat u_{i+1}^\prime=a(\alpha)$
with $\alpha\in H^{*-1}(\cE_{i+1};\tV )$. It follows that
$\kappa_{i}^{i+2,*} \tilde  u_{i+2}-\hat u_i=a(\kappa_i^*\alpha)$.
 By \ref{cikdcdc}
we can find a class
$\tilde \alpha\in H^{*-1}(\cE_{i+2};\tV )$
such that
$\kappa^{i+1,*}_{i}(\tilde\alpha)=\kappa^{*}_{i}(\alpha)$.
We set $\hat u_{i+1}:= \hat u_{i+1}^\prime+a(\alpha)-a(\kappa_{i+1}^*\tilde\alpha)$
and  $\hat u_{i+2}^\prime:=\tilde u_{i+2}-a(\tilde \alpha)$.
Then we have 
$\kappa_{i}^*\hat u_{i+1}=\hat u_{i}$ and
$\kappa_{i+1}^*(\hat u_{i+2}^\prime)=\hat u_{i+1}$.
\hB

\section{The natural transformation $\Phi$}
  
We consider a generalized cohomology theory $E$ represented by a spectrum $\bE$. We  consider two smooth extensions
$(\hat E,R,I,a)$ and $(\hat E^\prime,R^\prime,I^\prime,a^\prime)$ of $E$. Let us fix a degree
$k\in \Z$ and a classifying space $\bE_k$ for the homotopy functor $X\mapsto E^k(X)$, e.g. the $k$-th space of the spectrum $\bE$ if the latter is an $\Omega$-spectrum. In the present section we give a   construction of a natural transformation
$\Phi\colon \hat E^k\to \hat E^{\prime k}$ such that the following diagram
commutes:
$$\xymatrix{\Omega^{k-1}(M,\tV )\ar[r]^a\ar@{=}[d]&\hat E^k(M)\ar[d]^{\Phi}\ar[r]^{I}\ar@/^1cm/[rr]^R&E^k(M)\ar@{=}[d]&\Omega^k_{cl}(M,\tV )\ar@{=}[d]\\\Omega^{k-1}(M,\tV )\ar[r]^{a^\prime}& \hat E^{\prime k}(M)\ar[r]^{I^\prime}\ar@/_1cm/[rr]^{R^\prime}&E^k(M)&\Omega^k_{cl}(M,\tV )}
\ .$$
We make one of the following \textcolor{black}{two} assumptions:
\begin{ass}\label{uwidwqudqdqwidwqd}
\begin{enumerate}
\item\label{udwidw} $E^{k-1}=\pi_1(\bE_k)=0$, the abelian groups $E^m$ are finitely generated for all $m\le k$, and the smooth extensions are defined on the category of compact manifolds, or
\item\label{udwidw1} the smooth extensions are defined on the category of all manifolds and
 the abelian group $E^m$ is countably generated for all $m\le k$.
\end{enumerate}
\end{ass}
Note that $\pi_i(\bE_k)\cong E^{k-i}$. 
We choose an approximation $(\cE_i,x_i,\kappa_i)$ of $\bE_k$ by smooth manifolds as in Propositions \ref{duwidwdwd} or \ref{duwidwdwd1}. In case \ref{udwidw} we can assume that the $\cE_i$ are compact.

Let $u\in E^k(\bE_k)$ be the tautological class represented by the identity $\id\in [\bE_k,\bE_k]$.
We choose a  family of closed forms  $\omega_i\in \Omega^k_{cl}(\cE_i,\tV )$ as in Proposition \ref{fudif}
such that $\Rham(\omega_i)=\ch(x_i^*u)$ and $\kappa_i^*\omega_{i+1}=\omega_i$ for all $i\ge 0$. Then by Proposition \ref{89fwefwefewfwe} we choose families of classes $\hat u_i\in \hat E^k(\cE_i)$, $\hat u_i^\prime\in \hat E^{\prime k}(\cE_i)$ such that $R(\hat u_i)=R^\prime(\hat u_i^\prime)=\omega_i$
and $I(\hat u_i)=
I^\prime(\hat u_i^\prime)=x_i^*u$ and $\kappa_i^*\hat u_{i+1}=\hat u_i$, $\kappa_i^*\hat u_{i+1}^\prime=\hat u_i^\prime$ for all $i\ge 0$.

Now we can define the transformation $\Phi$, which  may depend on the choices made above.
Let $M$ be a (compact in case \ref{udwidw}) manifold  and $\hat v\in \hat E^k(M)$.
Note that $I(\hat v)=f^*u$ for some $f\in [M,\bE_k]$. By Property \ref{la0}  
there exists an $i\in \nat$ and a smooth map
$f_i\colon M\to \cE_i$ such that $f_i^*x_i^*u=I(\hat v)$.
Therefore there exists a unique $\alpha \in \Omega^{k-1}(M,\tV )/\im(\ch)$
such that $a(\alpha)+f_i^*(\hat u_i)=\hat v$.
We set
$$\Phi(\hat v):=a^\prime(\alpha)+f_i^*(\hat u_i^\prime)\ .$$

\begin{lem}
$\Phi$ is well-defined.
\end{lem}
\proof
The only choice involved is the map $f_i$.
We can increase the index $i$ to $j$ without changing $\Phi^k(\hat v)$
by replacing $f_i$ by $f_j:=\kappa^{j}_i\circ f_i$.
Given $f_i$ and $f_{i^\prime}^\prime$ then by property \ref{duwidwdwd}.\ref{la0} (or  \ref{duwidwdwd1}.\ref{la0}, resp.) there exists
$j\ge \max\{i,i^\prime\}$ such that $\kappa^j_i\circ f_i$ and $\kappa_{i^\prime}^j\circ f_{i^\prime}^\prime$ are 
homotopic.

Thus let us assume that we have $f_j$ and $f_j^\prime$ which are homotopic.
Let $H\colon I\times M\to \cE_j$ denote the homotopy.
We define
$\beta:=\int_{I\times M/M} H^*\omega_j\in \Omega_{cl}^{k-1}(M,\tV )/\im(\ch)$.
Then by the homotopy formula (\ref{udqwdqwdqw1}) we have 
$f^{\prime,*}_j(\hat u_j)=f^*_j(\hat u_j)+a(\beta)$, but also
$f^{\prime,*}_j(\hat u_j^\prime)=f^*_j(\hat u_j^\prime)+a^\prime(\beta)$.
If $\alpha^\prime$ and $\Phi^{\prime}(\hat v)$ denote the   result  for  $\alpha$  and $\Phi(\hat v)$  for the choice $f_j^\prime$, then we have
$\alpha^\prime=\alpha-\beta$. But this implies 
$$\Phi^{\prime}(\hat v)=a^\prime(\alpha^\prime)+f^{\prime,*}_j(\hat u^\prime_j)=a^\prime(\alpha)-a^\prime(\beta)+f^{\prime,*}_j(\hat u_j^\prime)=
a^\prime(\alpha)+f_i^*(\hat u_i^\prime)=\Phi(\hat v)\ .$$\hB

 \begin{lem}
We have by construction $R^\prime\circ \Phi=R$,
$I^\prime\circ \Phi=I$, and $\Phi\circ a=a^\prime$.
 \end{lem}
 \proof
Straightforward verifications.
\hB

\begin{lem}
$\Phi$ is natural. 
\end{lem}
\proof
For a moment we write $\Phi_M$, $\Phi_N$ in order to indicate a possible dependence on the underlying manifold. 
Let $g\colon N\to M$ be a smooth map between  manifolds. Given
$\hat v\in \hat E^k(M)$ we must show that
$g^*\Phi_M(\hat v)=\Phi_N(g^*\hat v)$.
Note that $\Phi_M(\hat v)=a^\prime(\alpha)+f_i^*(\hat u_i^\prime)$,
$g^*\hat v=a(g^*\alpha)+g^*f_i^*(\hat u_i)$, and therefore
$\Phi_N(g^*\hat v)=a^\prime(g^*(\alpha))+g^*f_i^*(\hat u_i^\prime)=g^*\Phi_M(\hat v)$.
\hB 

Note that $\Phi$ as defined above is a natural transformation of set-valued functors.
In general it does not preserve the group structures. 
The deviation from additivity is a priory a natural transformation
$$\tilde B\colon \hat E^k\times \hat E^k\to \hat E^{\prime k}$$ 
such that
\begin{equation}\label{ezwu}
\Phi(\hat v+\hat w)=\Phi(\hat v)+\Phi(\hat w)+\tilde B(\hat v,\hat w)\ . 
\end{equation}
which satistifies the cocycle condition and symmetry
$$\tilde B(\hat u,\hat v+\hat w)+\tilde B(\hat v,\hat w)=\tilde B(\hat u,\hat
v)+\tilde B(\hat u+\hat v,\hat w)\ , \quad \tilde B(\hat u,\hat v)=\tilde B(\hat v,\hat u)\ .$$
Because of the identities 
 \begin{eqnarray}
 0&=&\Phi(\hat v+a(\alpha))-a^\prime(\alpha)-\Phi(\hat v)=\tilde B(\hat v,a(\alpha))\nonumber\\
 0&=&R^\prime(\Phi(\hat v+\hat w)-\Phi(\hat v)-\Phi(\hat w))=R^\prime(\tilde
 B(\hat v,\hat w))\label{ide1}\\ 0&=&I^\prime(\Phi(\hat v+\hat w)-\Phi(\hat
 v)-\Phi(\hat w))=I^\prime(\tilde B(\hat v,\hat w))\nonumber
 \end{eqnarray}
it factors over a natural transformation 
\begin{equation}\label{uewiqeeqw}
B\colon E^k(M)\times E^k(M)\to H^{k-1}(M; \tV )/\im(\ch)\ .
\end{equation}

%
%
\begin{ddd}\label{uzdiqwdqwdwqd666}
We call the cohomology \textcolor{black}{theory} $E$ \textbf{rationally even}, if $E^m\otimes_\Z\Q=0$ for all odd $m\in \Z$.
\end{ddd}

\begin{theorem}\label{hjsasas}
Let $k\in \Z$ be even.  If $E^*$ is rationally even and one of \ref{uwidwqudqdqwidwqd}.\ref{udwidw} or \ref{uwidwqudqdqwidwqd}.\ref{udwidw1} is satisfied,
then the transformation
$\Phi\colon \hat E^k\to \hat E^{k,\prime}$ is additive. 
\end{theorem}
\proof
The family $(\cE_i\times \cE_i,\kappa_i\times \kappa_i)_{i\ge 0}$ of manifolds gives rise to a system of abelian groups $(H^{k-1}(\cE_i\times \cE_i;\tV ),(\kappa_i\times \kappa_i)^*)_{i\ge 0}$ indexed by $\nat^{op}$.
The natural transformation $B$ induces a class
$$\hat B\in \lim_i ( H^{k-1}(\cE_i\times \cE_i;\tV )/\im(\ch))\ .$$
In detail, the $i$th component is the class
$$\hat B_i:=B(\pr_1^*x_i^*u,\pr_2^*x_i^*u)\in H^{k-1}(\cE_i\times \cE_i;\tV )/\im(\ch)\ .$$

\begin{prop} 
We have $\hat B=0$.
\end{prop}
\proof
We show this result by showing that $\lim ( H^{k-1}(\cE_i\times \cE_i;\tV )/\im(\ch))=0$.
 We first show a refinement of the Mittag-Leffler condition. We start with the following  general fact.
\begin{lem}\label{ghdgfd}
Let $X$ be a  $CW$-complex such that $\pi_{2i+1}(X)\otimes_\Z \Q=0$\footnote{If $X$ is not connected, then we require this for all its components.} for $i=0,\dots,n$.
Then $H_{2i+1}(X;\Q)=0$ for $i=0,\dots,n$.
\end{lem}
\proof
We first show that $H^{2i+1}(X;\Q)=0$ for $i=0,\dots,n$.
We assume the contrary.
Let $k\in \{0,\dots,n\}$ be the smallest number
such that $H^{2k+1}(X;\Q)\not=0$. We have 
the Postnikov tower
 $$X\langle 2k+2\rangle \to X\langle 2k+1\rangle \to \dots\to X\langle 2\rangle \to X\langle 1\rangle \to X\ .$$
The fibre of $X\langle 2j\rangle \to X\langle 2j-1\rangle $ is equivalent to an Eilenberg-MacLane space  $K(\pi_{2j-1}(X),2j-2)$
which is rationally acyclic since $\pi_{2j-1}(X)\otimes_\Z \Q=0$.
Therefore
$X\langle 2j\rangle \to X\langle 2j-1\rangle $ induces a rational cohomology equivalence. 
The fibre of $X\langle 2j+1\rangle \to X\langle 2j\rangle $ is  an Eilenberg-MacLane space $K(\pi_{2j}(X),2j-1)$.
With the Serre spectral sequence
$$H^*(X\langle 2j\rangle ;H^*(K(\pi_{2j}(X),2j-1);\Q))\Longrightarrow H^*(X\langle 2j+1\rangle ;\Q)$$
and 
$$H^l(K(\pi_{2j}(X),2j-1),\Q)=\left\{\begin{array}{cc} 0&l\not\in \{0,2j-1\}  \\ \Q&l=0,2j-1\end{array}
\right.$$  we conclude that
$$H^{2k+1}(X;\Q)\to H^{2k+1}(X\langle 1\rangle ;\Q)\to \dots\to H^{2k+1}(X\langle  2k+2\rangle ;\Q)=0$$ is injective.  This is a contradiction. 

It now follows by duality that
$H_{2i+1}(X;\Q)\cong 0$ for $i\in \{0,\dots ,n\}$.
\hB

\begin{lem}\label{zuwqeqe}
For every $i\in \nat$ and  all $r\ge i+1$ we have 
\begin{equation*}
(\kappa^{i+r}_{i}\times \kappa_i^{i+r})^* H^{k-1}(\cE_{i+r}\times \cE_{i+r};\tV )=0\ .
\end{equation*}
\end{lem}
\proof
For $r\ge i+1$ the map $x_{i+r}\colon \cE_{i+r}\to \bE_k$ is $2i+1$-connected so that
$\pi_{l}(\cE_{i+r})\cong \pi_l( \bE_k)\cong E^{k-l}$ for all
$l\le 2i$. In particular, $\pi_{2l-1}(\cE_{i+r})\otimes_\Z\Q=0$ and hence
$\pi_{2l-1}(\cE_{i+r}\times \cE_{i+r})\otimes_\Z\Q=0$ for all $l\in \nat$ such
that $2l-1\le 2i$.   
It follows from Lemma
\ref{ghdgfd} that $H^{2l-1}(\cE_{i+r}\times \cE_{i+r};\Q)=0$
for all $l\in \nat$ such that $2l-1\le 2 i$. Since $\cE_i\times \cE_i$ is homotopy equivalent to a $2i$-dimensional complex
and $\tV^j=0$ for odd $j$
 we now have
\begin{eqnarray*}\lefteqn{
(\kappa^{i+r}_{i}\times \kappa_i^{i+r})^* H^{k-1}(\cE_{i+r}\times \cE_{i+r};\tV )}&&\\&\cong &
(\kappa^{i+r}_{i}\times \kappa_i^{i+r})^* \bigoplus_{2l-1\le 2i} H^{2l-1}(\cE_{i+r}\times \cE_{i+r};\tV^{k-2l})\\
&=&0\ .
\end{eqnarray*}
\hB 

We now consider the system of exact sequences
$$0\to \ch(E^{k-1}(\cE_i\times \cE_i))\to H^{k-1}(\cE_i\times \cE_i;\tV )\to\frac{H^{k-1}(\cE_i\times \cE_i;\tV )}{ \ch(E^{k-1}(\cE_i\times \cE_i))}\to 0$$
indexed by $\nat^{op}$.
We apply $\lim$ and get the exact sequence 
$$\lim_i (H^{k-1}(\cE_i\times \cE_i;\tV )) \to \lim_i \left(\frac{H^{k-1}(\cE_i\times \cE_i;\tV )}{ \ch(E^{k-1}(\cE_i\times \cE_i))}\right)\to {\lim_i}^1 (\ch(E^{k-1}(\cE_i\times \cE_i)))\ .$$
Now  we have $\lim_i ( H^{k-1}(\cE_i\times \cE_i;\tV ))=0$
because of Lemma \ref{zuwqeqe}. The same lemma implies that
the subsystem  $\ch(E^{k-1}(\cE_i\times \cE_i))\subseteq  H^{k-1}(\cE_i\times
\cE_i;\tV )$ satisfies the Mittag-Leffler condition so that $\lim_i^1
(\ch(E^{k-1}(\cE_i\times \cE_i)))=0$.
This implies  $$\lim_i  \left(\frac{H^{k-1}(\cE_i\times \cE_i;\tV )}{ \ch(E^{k-1}(\cE_i\times \cE_i))}\right)=0$$ as required.
\hB

We now show Theorem \ref{hjsasas}.
Let ${ \hat v }, { \hat w }\in \hat E^k(M)$.
We choose $i$ sufficiently large such that there exist
$f_{ \hat v },f_{ \hat w }\colon M\to \cE_i$ with $f_{ \hat v }^*(x_i^*u)=I({ \hat v })$ and $f_{ \hat w }^*(x_i^*u)=I({ \hat w })$.
Let $j\ge i$ be such that there exists $\mu\colon \cE_i\times \cE_i\to \cE_j$ with
$\mu^*(x_j^*u)=\pr_0^*x_i^*u+\pr_1^*x_i^*u$. We choose $\alpha\in \Omega^{k-1}(\cE_i\times \cE_i,\tV )$ such that
$a(\alpha)+\mu^*(\hat u_j)=\pr_0^*\hat u_i+\pr_1^*\hat u_i$.

We further can choose $f_{{ \hat v }+{ \hat w }}\colon M\to \cE_j$
as the composition
$f_{{ \hat v }+{ \hat w }}=\mu\circ (f_{ \hat v },f_{ \hat w })$ so that $f_{{ \hat v }+{ \hat w }}^*(x_j^*u)=I({ \hat v }+{ \hat w })$.
 
 We now choose $\alpha_{ \hat v },\alpha_{ \hat w },\alpha_{{ \hat v }+{ \hat w }}\in \Omega^{k-1}(M,\tV )$ such that
$a(\alpha_{ \hat v })+f_{ \hat v }^*(\hat u_i)={ \hat v }$,  $a(\alpha_{ \hat w })+f_{ \hat w }^*(\hat u_i)={ \hat w }$, and
$a(\alpha_{{ \hat v }+{ \hat w }})+f_{{ \hat v }+{ \hat w }}^*(\hat u_j)={ \hat v }+{ \hat w }$.

Then we have
$\Phi({ \hat v })=a^\prime(\alpha_{ \hat v })+f_{ \hat v }^*( \hat  u_i^\prime)$, 
$\Phi({ \hat w })=a^\prime(\alpha_{ \hat w })+f_{ \hat w }^*(\hat u_i^\prime)$, and
$\Phi({ \hat v }+{ \hat w })=a^\prime(\alpha_{{ \hat v }+{ \hat w }})+f_{{ \hat v }+{ \hat w }}^*(\hat u_j^\prime)$.
We now calculate using $\Phi(\hat u_i)=\hat u_i^\prime$ 
and 
$$0=B(\pr_0^*x_i^*u,\pr_1^*x_i^*u)=\Phi(\pr_0^*\hat u_i+\pr_1^*\hat u_i)-\Phi(\pr_0^*\hat u_i)-\Phi(\pr_1^*\hat u_i)$$  at the marked equality
\begin{eqnarray*}
\lefteqn{\Phi({ \hat v }+{ \hat w })-\Phi({ \hat v })-\Phi({ \hat w })}&&\\&=&a^\prime(\alpha_{{ \hat v }+{ \hat w }})+f_{{ \hat v }+{ \hat w }}^*(\hat u_j^\prime)-
a^\prime(\alpha_{ \hat v })-f_{ \hat v }^*( \hat  u_i^\prime)-a^\prime(\alpha_{ \hat w })-f_{ \hat w }^*(\hat u_i^\prime)\\
&=&a^\prime(\alpha_{{ \hat v }+{ \hat w }})+(f_{ \hat v },f_{ \hat w })^*\mu^*(\hat u_j^\prime)-
a^\prime(\alpha_{ \hat v })-f_{ \hat v }^*( \hat  u_i^\prime)-a^\prime(\alpha_{ \hat w })-f_{ \hat w }^*(\hat u_i^\prime)\\&=&
a^\prime(\alpha_{{ \hat v }+{ \hat w }}-\alpha_{ \hat v }-\alpha_{ \hat w })+(f_{ \hat v },f_{ \hat w })^*\left(\Phi(\pr_0^*\hat u_i+\pr_1^*\hat u_i)-a^\prime(\alpha)\right)\\&&-
 f_{ \hat v }^*( \hat  u_i^\prime)-f_{ \hat w }^*(\hat u_i^\prime)\\
&\stackrel{!}{=}&
a^\prime(\alpha_{{ \hat v }+{ \hat w }}-\alpha_{ \hat v }-\alpha_{ \hat w })+(f_{ \hat v },f_{ \hat w })^*\left(\Phi(\pr_0^*\hat u_i)+\Phi(\pr_1^*\hat u_i)-a^\prime(\alpha)\right)\\&&-
 f_{ \hat v }^*( \hat  u_i^\prime)-f_{ \hat w }^*(\hat u_i^\prime)\\
&=&a^\prime(\alpha_{{ \hat v }+{ \hat w }}-\alpha_{ \hat v }-\alpha_{ \hat w })+f_{ \hat v }^*\hat u_i^\prime+f_{ \hat w }^*\hat u_i^\prime -(f_{ \hat v },f_{ \hat w })^* a^\prime(\alpha)\\&&-
 f_{ \hat v }^*( \hat  u_i^\prime)-f_{ \hat w }^*(\hat u_i^\prime)\\
&=&a^\prime(\alpha_{{ \hat v }+{ \hat w }}-\alpha_{ \hat v }-\alpha_{ \hat w }-(f_{ \hat v },f_{ \hat w })^*\alpha)\ .
\end{eqnarray*}
Doing the same calculation starting with $0=({ \hat v }+{ \hat w })-{ \hat v }-{ \hat w }$ (leave out the symbols
$\Phi$ and ${}^\prime$)
we get $0=a(\alpha_{{ \hat v }+{ \hat w }}-\alpha_{ \hat v }-\alpha_{ \hat w }-(f_{ \hat v },f_{ \hat w })^*\alpha)$. 
Since $\ker(a)=\ker(a^\prime)$ we conclude that $$\Phi({ \hat v }+{ \hat w })-\Phi({ \hat v })-\Phi({ \hat w })=0\ .$$
\hB 

\begin{theorem}\label{main1}
Let $E$ be a rationally even generalized  cohomology theory  which is  represented by a spectrum $\bE$.  Let
$(\hat E,R,I,a,\int)$ and $(\hat E^\prime,R^\prime,I^\prime,a^\prime,\int^\prime)$ be  two smooth extensions with integration. We assume that either the smooth extensions are defined on the category of all smooth manifolds and the groups $E^m$ are countably generated for all $m\in \Z$, or that $E^m=0$ for all odd $m\in \Z$ and $E^m$ is finitely generated for even $m\in \Z$.
Then there is a unique natural isomorphism
$$\Phi\colon \hat E\to \hat E^\prime$$ of smooth extensions with integration.
\end{theorem}
\proof
Let us first show the existence of a natural transformation. 
We let $\Phi^k\colon \hat E^k\to \hat E^{\prime k}$ denote the component in degree $k$.
Let $\Phi^{2k}$ be the transformation obtained in Theorem \ref{hjsasas}. We extend $\Phi$
to odd degrees using integration.

Let $i\colon M\to S^1\times M$ be the embedding  induced by a point in $S^1$ and $p\colon S^1\times M\to M$ be the projection. Because of $p\circ i=\id$ we have a splitting $$E^*(S^1\times M)\cong p^* E^*(M)\oplus \ker(i^*)\ .$$ Let $q\colon S^1\times M\to \Sigma M_+$ be the projection onto the suspension. It induces an identification
$\ker(i^*)\cong q^*\tilde E^*(\Sigma M_+)\cong E^{*-1}(M)$ using the suspension isomorphism.
We let $\sigma\colon E^{*-1}(M)\to E^*(S^1\times M)$ be the inclusion of this summand.

Assume that we already have constructed $\Phi^{k}$. This is in particular the case for even $k$.
Let
$\hat x\in \hat E^{k-1}(M)$ be given.
Then we choose $\tilde x\in  \hat E^{k}(S^1\times M)$ such that
\begin{enumerate}
\item \label{uierwer0} $\int(\tilde x)=\hat x$,
\item\label{uierwer1} $R(\tilde x)=dt\wedge \pr^*R(\hat x)$,
\item\label{uierwer2} $I(\tilde x)=\sigma(I(\hat x))$.
\end{enumerate}
The following procedure shows that this choice can be made.
We first choose a lift $\tilde x\in \hat E^{k}(S^1\times M)$ of $\sigma( I(\hat x))\in E^{k}(S^1\times M)$ so that \ref{uierwer2} is satisfied. Then we add $a(\alpha)$
for a suitable $\alpha\in \Omega^{k-1}(S^1\times M,\tV )$ in order to satisfy \ref{uierwer1}.
Then $\hat x-\int(\tilde x)\in a(H^{k-1}(M;\tV ))$. We can kill this difference by modifying
$\alpha$ suitably. Here we use that $\int(a(\alpha))=a(\int(\alpha))$ and
$\int\colon \Omega^*(S^1\times M,\tV )\to \Omega^{*-1}(M,\tV )$ is surjective.

 We  now define
$$\Phi^{k-1}(\hat x):=\int^\prime(\Phi^{k}(\tilde x))\ .$$
Let us check  that this is well-defined.
Note that another choice $\tilde x^\prime$ satisfies
$\tilde x^\prime-\tilde x=a(\alpha)$ with $\int(a(\alpha))=0$. Furthermore,
$$\int^\prime(\Phi^{k}(\tilde x^\prime))=\int^\prime(\Phi^{k}(\tilde x+a(\alpha)))=\int^\prime(\Phi^{k}(\tilde x)+a(\alpha))=\int^\prime( \Phi^{k}(\tilde x))\ .
$$
Naturality of $\Phi^{k-1}$ follows from naturality of the integration maps.
Indeed, let $f\colon N\to M$ be a smooth map and $\hat y:=f^*\hat x$. Then
we can choose $\tilde y:=(\id_{S^1}\times f)^*\tilde x$. We get
$$\Phi^{k-1}(\hat y)=\int^\prime(\Phi^{k}(f^*\tilde x))=\int^\prime
(\id_{S^1}\times f)^*\Phi^{k}(\tilde x)= f^* \int^\prime (\Phi^{k}(\tilde
x))=f^*\Phi^{k-1}(\hat x)\ .$$

Let us now discuss uniqueness of $\Phi^{2k}$.
Assume that $\Psi\colon \hat E^{2k}\to \hat E^{\prime 2k}$ is a second natural transformation
of group-valued functors which is compatible with the transformations $R,I,a$ in the sense that
$$R^\prime \circ \Psi=R\ ,\quad I^\prime\circ \Psi=I\ , \quad a^\prime = \Psi\circ a\ .$$
Then we consider the difference
$\Delta:=\Phi^{2k}-\Psi$. Compatibility with $R$ shows that $\Delta$ takes values in
$\hat E^{\prime 2k}_{flat}$ (see Definition \ref{uiqwdqwdwqd}). Compatibility with $I$ in addition shows that
$\Delta$ takes values in the subfunctor $a^\prime(H^{2k-1}(\dots;\tV )/\im(\ch))\subseteq \hat E^{\prime 2k}$. Finally, compatibility with $a$ implies that $\Delta$ descends to a natural transformation
$$\Delta\colon E^{2k}\to H^{2k-1}(\dots;\tV )/\im(\ch)\ .$$
We get an element
$$(\Delta(x_i^*u))\in \lim_i( H^{2k-1}\left(\cE_i;\tV )/\im(\ch))\right)\ .$$
The same argument as for Lemma \ref{zuwqeqe} shows 
that the target group vanishes so that $\Delta(x_i^*u)=0$ for all $i\in \nat$. But this implies that $\Delta=0$. Indeed, if $M$ is some manifold and $x\in E^{2k}(M)$, then there exists an $i\in \nat$ and $f\colon M\to \cE_i$ such that
$x=f^*x_i^*u$. We get $\Delta(x)=f^*\Delta(x_i^*u)=0$. 
 
Recall that we have defined $\Phi^{2k}$ using Theorem \ref{hjsasas}. Then we
have extended the construction to odd degrees such that $\Phi^{2k-1}$ is
chracterized by
\begin{equation}\label{duidqwdqwd}
\int^\prime \circ\:\: \Phi^{2k}=\Phi^{2k-1}\circ \int\ .
\end{equation}
We could use the construction above in order to construct another transformation
$\Psi\colon \hat E^{2k}\to \hat E^{\prime  2k}$ starting from $\Phi^{2k+1}$ so that
$$\int^\prime\circ \:\: \Phi^{2k+1}=\Psi\circ \int\ .$$
By the uniqueness result we see that $\Psi=\Phi^{2k}$.

Therefore we have constructed a natural transformation of smooth extensions with integration. As such it is unique on the even part.
Since the integration $\int\colon \hat E^{2k}(S^1\times M)\to \hat E^{2k-1}(M)$ is surjective, the compatibility (\ref{duidqwdqwd}) implies uniqueness on the odd part, too. 
\hB

\section{Multiplicative structures}
In this section
we assume that $E$ is a multiplicative cohomology theory.
Let $(\hat E,R,I,a)$ be a smooth extension of $E$.
We fix the unique $e\in \ker(i^*)\subseteq E^1(S^1)$ such that $\int_{E}(e)=1\in E^0$, 
the unit of the ring $E^*$.
We let $\omega_e\in \Omega^1_{cl}(S^1;\tV^0 )$ be the unique  rotationally invariant closed form
such that $\Rham(\omega_e) =\ch(e)$. In coordinates $\omega_e=dt\otimes 1$.
Note that $\int_{S^1} \omega_e=1$.
Then we choose a lift $\hat e\in \hat E^1(S^1)$ such that $I(\hat e)=e$ and $R(\hat e)=\omega_e$.
Note that $\hat e$ is determined uniquely up to elements of the form
$a(\alpha)$ with $\alpha\in H^{0}(S^1;\tV )/\im(\ch)$.
We can assume that the representative $\alpha$ is constant.
We see that $\hat e$ is determined uniquely up to a choice
in $\tV^0/E^0\oplus \tV^{-1}/E^{-1}$.

We now want to modify the class $\hat e$ such that $q^*\hat e=-\hat e$, where $q\colon S^1\to S^1$ is  given by
$q(z)=\bar z$. Apriori
$R(q^*\hat e+\hat e)=0$ and $I(q^*\hat e+\hat e)=0$. Therefore
$q^*\hat e+\hat e=a(\rho)$ for $\rho\in H^0(S^1;\tV )$.
We write $\rho=\rho_0+\rho_1\in H^0(S^1;\tV^0)\oplus H^1(S^1;\tV^{-1})$.
Since $a(q^*\rho-\rho)=0$, $q^*\rho_0=\rho_0$ and $q^*\rho_1=-\rho_1$ we conclude that $2\Rham(\rho_1)=
\ch(r)$ for some $r\in E^0(S^1)$.
Then
$$q^*(\hat e-a(\frac{1}{2}\rho))+ \hat e-a(\frac{1}{2}\rho)=a(\rho-\frac{1}{2}\rho-\frac{1}{2}q^*\rho)=a(\frac{1}{2}(\rho-q^*\rho))=a(\rho_1)=a(\frac{1}{2}\ch(r))\ .$$ If we replace $\hat e$ by $\hat e-a(\frac{1}{2}\rho)$, then the new $\hat e$ satisfies
$q^*\hat e+\hat e=a(\rho_1)$.
The $2$-torsion class
$[\Rham(\rho_1)] \in \tV^{-1}/E^{-1}$ is independent of
the choices.
Indeed, if we change
$\hat e$ by $\hat e+a(\alpha)$ for $\alpha\in H^0(S^1;\tV)$, then we can take the same
$\rho_1$.

Recall the inclusion $i\colon M\to S^1\times M$ and the projection $p\colon S^1\times M\to M$.
We finally consider the condition $i^*\hat e=0$. Of course,  we have $i^*\hat e=a(\theta)$ for some $\theta\in H^0(*;\tV)\cong \tV^0$. If we replace $\hat e$ by $\hat e-p^*\theta$, then we get $i^*\hat e=0$
retaining all the other conditions.  The ramaining choice for $\hat e$ is $\tV^{-1}/E^{-1}$.

\begin{ddd}
The class $\bo_{\hat E}:=[\Rham(\rho_1)]\in \tV^{-1}/E^{-1}$ is called the \textbf{obstruction class}.
\end{ddd}
 The obstruction class vanishes exactly if we can choose  $\hat e\in \hat E^1(S^1)$ such that
$$R(\hat e)=dt\otimes 1\ , \quad I(\hat e)=e\ ,\quad i^*\hat e=0\ ,\quad q^*\hat e=-\hat e\ .$$
In this case $\hat e$ is unique up to an element in $\tV^{-1}/E^{-1}$.
The obstruction vanishes e.g. if $E^{-1}$ is torsion, and in this case $\hat e$ is unique.
We do not have any example of a smooth extension with non-trivial obstruction class.

\begin{prop}\label{uzqwiwqdqd} If  $(\hat E,R,I,a)$ is a multiplicative smooth extension with vanishing obstruction class $\bo_{\hat E}$, then it has an integration. 
\end{prop}
 \proof

Using the class $e$ we can make the decomposition
$$E^{*+1}(S^1\times M)\cong \im(p^*)\oplus ¸\ker(i^*)\cong E^{*+1}(M)\oplus E^*(M)$$
more explicit. Namely, we can write the class
$x\in E^{*+1}(S^1\times M)$ uniqely in the form
$x=p^*u\oplus e\times y$ with $u\in E^{*+1}(M)$ and $y\in E^*(M)$.
Note that $y=\int x$.

We now use the decomposition
 $$\hat E^{*+1}(S^1\times M)=\im(p^*)\oplus \ker(i^*)$$
in order to
define an integration  which factors as
$$\int \colon \hat E^{*+1}(S^1\times M)\stackrel{\pr}{\to} \ker(i^*) \to \hat E^*(M)\ .$$ 
It  obviously satisfies the second condition \ref{ddd4}.\ref{ddd41}
$$\int\circ\:\: p^*=0\ .$$ 
Let $\hat x\in \ker(i^*)$. Then we can write
$I(\hat x)=e\times y$ for a unique $y\in E^*(M)$.
We choose a smooth lift $\hat y\in \hat E^*(M)$ and a form
$\rho\in \Omega^{*-1}(S^1\times M,\tV)$ such that
$\hat x=\hat e\times \hat y +a(\rho)$. 
Then we define
$$\int \hat x:= \hat y+a(\int  \rho)\ .$$
Let us show that $\int$ is well-defined.
If we choose another smooth lift $\hat y^\prime$, then
$\hat y^\prime=\hat y+a(\theta)$ for some form $\theta\in \Omega^{*-1}(M,\tV)$.
Since $\hat e\times \hat y^\prime=\hat e\times \hat y+\hat e\times a(\theta)=\hat e\times \hat y+a(\omega_e\times \theta)$ we can choose $\rho^\prime=\rho-\omega_e\times \theta$ and  our construction produces
$$\hat y^\prime + a(\int \rho^\prime)=\hat y+a(\theta)+a(\int  \rho)-a(\int  \omega_e\wedge\theta)=\hat y+a(\int \rho)\\ .$$
Next we show that the construction is independent of the choice of $\rho$.
If we choose $\rho^\prime$, then
$\rho^\prime=\rho+\alpha$ such that $\alpha\in \Omega^{*-1}(S^1\times M,\tV)$ is closed and $\Rham(\alpha)=\ch(u)$ for some $u\in E^{*}(S^1\times M)$.
But then our construction produces
$$\hat y+a(\int  \rho^\prime)=\hat y+a(\int  \rho)+a(\int  \ch(u))=\hat y+a(\int  \rho)+a(\ch(\int u))=\hat y+a(\int  \rho)\ .$$
 
We have
$$I(\int \hat x)=I(\hat y+a(\int \rho))=I(\hat y)=y=\int I(\hat x)\ .$$
Furthermore
$$R(\int \hat x)=R(\hat y+a(\int\rho))=\int(\omega_e\times R(\hat y))+d\int\rho=\int(\omega_e\times R(\hat y)+d\rho)=\int R(\hat x)\ .$$

Next we check that the integration is linear.
Let $\hat x_i\in  \ker(i^*)$, $i=0,1$. We choose $y_i\in E^*(M)$ such that $I(x_i)=e\times y_i$
and smooth lifts $\hat y_i\in \hat E^*(M)$. Then we find $\rho_i\in  \Omega^{*-1}(S^1\times M,\tV)$ such that
$\hat x_i=\hat e\times \hat y_i+a(\rho_i)$ for $i=0,1$.  Note that
$\hat x_0+\hat x_1=\hat e\times(\hat y_0+\hat y_1)+a(\rho_0+\rho_1)$. It follows
$$\int  (\hat x_0+\hat x_1)=\hat y_0+\hat y_1+a(\int(\rho_0+\rho_1))=\int (\hat x_0)+\int (\hat x_1)\ .$$

If $\hat x=\hat e\times \hat y+a(\rho)$, then
$$(q\times \id_M)^*\hat x=q^*\hat e\times \hat y+a((q\times \id_M)^*\rho)=\hat e\times (-\hat y)+a((q\times \id_M)^*\rho)\ .$$
We get
$$\int((q\times \id_M)^*\hat x)=-\hat y+a(\int (q\times \id_M)^*\rho) =-\hat y-a(\int \rho )=-\int \hat x \ .$$


Let us finally show that the integration is natural.
Let $f\colon N\to M$ be a smooth map. Then we can write
$f^*\hat x=\hat e\times f^*\hat y+(\id_{S^1}\times f)^* \rho$.
It follows that
$$\int (f^*\hat x)=f^*\hat y+a(f^*\int  \rho)=f^*\int (\hat x)\ .$$
\hB 

If we change $\hat e$ to $\hat e^\prime=\hat e+a(u)$, $u\in H^1(S^1;\tV^{-1})/\im(\ch)\cap H^1(S^1;\tV^{-1})\cong \tV^{-1}/E^{-1}$, then we have
$$\hat x=\hat e^\prime\times\hat y+a(\rho^\prime)=\hat e\times \hat y+a(u\times R(\hat y)+\rho^\prime)\ .$$
For the corresponding integration we get, using $\Rham(R(\hat
y))=\ch(\int(I(\hat x)))$,
$$\int\hat x=\int^\prime\hat x+a(\int (u\times R(\hat y)))=\int^\prime\hat
x+a(u\times  \ch(\int I(\hat x) ))\ .$$
 If $\tV^{-1}\not=0$, then we indeed may change the integration by modifying the choice of $\hat e$.
If $E$ is rationally even, or more generally, if only $E^{-1}$ is a torsion group, then of course $\tV^{-1}=0$.
\begin{kor}\label{uifqfqfqf}
If $(\hat E,R,I,a)$ is a multiplicative extension of a  generalized cohomology theory such that $E^{-1}$ is a torsion group, then there is a canonical choice of an integration.
\end{kor}

We can now apply Theorem \ref{main1}. 
\begin{kor}\label{main2}
Let $(\hat E,R,I,a)$ and  $(\hat E^\prime,R^\prime,I^\prime,a^\prime)$ be two multiplicative extensions of a rationally even generalized cohomology theory. We assume that either both extensions are defined on the category of all smooth manifolds and the groups $E^m$ are countably generated for all $m\in \Z$, or they are defined on the category of compact manifolds , $E^m=0$ for all odd $m\in \Z$ and $E^m$ is finitely generated for even $m\in \Z$. Then there is a unique natural isomorphism between these  smooth extensions preserving the canonical integration. 
\end{kor}

\begin{lem}\label{dqedwqdwqd}
The integration defined above satisfies the projection formula
$$\int (\hat x\cup p^*\hat z)=(\int\hat x)\cup \hat z\ ,\quad \hat x\in \hat E^*(S^1\times M)\ ,\quad \hat z\in \hat E^*(M)\ .$$
\end{lem}
\proof
We write $\hat x=p^*\hat u+\hat e\times \hat y+  a(\rho)$.
Then by construction $\int \hat x=\hat y+a(\int\rho)$.
Furthermore
$$\hat x\cup p^*\hat z=p^*(\hat u\cup \hat z)+\hat e\times (\hat y\cup \hat z)+a(\rho\wedge R(p^*\hat z))\ .$$
It follows that
$$\int(\hat x\cup p^*\hat z)=\hat y\cup \hat z+a(\int\rho\wedge p^*R(\hat z))=\hat y\cup \hat z+a((\int\rho) \wedge R(\hat z))=(\int \hat x)\cup \hat z\ .$$
\hB

\begin{theorem}\label{main23}
The  unique natural transformation  of Corollary \ref{main2} is multiplicative. 
\end{theorem}
\proof
The argument is similar to the proof of Theorem \ref{hjsasas}. 
We transport the $\cup$-product of $\hat E^\prime$ to $\hat E$ using the natural transformation $\Phi\colon \hat E\to \hat E^\prime$. Thus we define
$$\cup^\prime\colon \hat E^{ev}(M)\otimes \hat E^{ev}(M)\to \hat E^{ev}(M)$$
by
$$\hat x\cup^\prime \hat y:=\Phi^{-1}(\Phi(\hat x)\cup \Phi(\hat y))\ .$$
The difference
$$\Delta(\hat x,\hat y):=\hat x\cup \hat y-\hat x\cup^\prime \hat y$$
is a natural transformation of set-valued functors
$$\hat E^{ev}\times \hat E^{ev}\to \hat E^{ev}\ .$$
Since both $\cup$-products are compatible with $R$ and $I$ we conclude that
$\Delta$ actually has values in the subfunctor $a(H^{odd}(\dots;\tV)/\im(\ch))$.
Furthermore, since both cup-products are compatible with $a$ the transformation factors as
$$\Delta\colon  E^{ev}\times E^{ev }\to a(H^{odd}(\dots;\tV)/\im(\ch))\ .$$
We approximate $\bE_k$ by a family of manifolds $\cE_{k,i}$ as in Proposition \ref{duwidwdwd}.
Then we consider compatible families
$\hat u_{k,i}\in \hat E^k(\cE_{k,i})$ as in Proposition \ref{89fwefwefewfwe}.
For even $k,l$  we get a family of elements
$$(\Delta(\hat u_{k,i},\hat u_{l,i}))_{i\ge 0}\in \lim\limits_{i} (H^{odd}(\cE_{k,i}\times \cE_{l,i})/\im(\ch))\ .$$
The analog of Lemma \ref{zuwqeqe} shows that this limit is trivial, so that
$\Delta(\hat u_{k,i},\hat u_{l,i})=0$ for all $i$.
As in the proof of Theorem \ref{main1} we conclude that this implies $\Delta=0$.

We now discuss multiplicativity in general.
Let $k$ be odd and $l$ be even, $\hat x\in \hat E^{k}(M)$ and $\hat y\in \hat E^l(M)$.
Then we have, using the projection formula \ref{dqedwqdwqd} and the compatibility of $\Phi$ with integration,
\begin{eqnarray*}
\Phi(\hat x\cup \hat y)&=&\Phi(\int(\hat e\times (\hat x\cup \hat y)))\\&=&\int (\Phi(\hat e\times (\hat x\cup \hat y)))\\&=&\int(\Phi((\hat e\times \hat x) \cup p^* \hat y))\\&=& 
\int(\Phi(\hat e\times \hat x)\cup \Phi(p^*\hat y))\\&=&\int(\Phi(\hat e\times \hat x) \cup p^*\Phi( \hat y))\\&=&
(\int\Phi(\hat e\times \hat x))\cup  \Phi(\hat y)\\
&=&\Phi(\int(\hat e\times \hat x))\cup \Phi(\hat y)\\
&=&\Phi(\hat x)\cup \Phi(\hat y) \ .
\end{eqnarray*}
Similarly, if $k$ and $l$ are odd, then again using the projection formula and the case just shown
\begin{eqnarray*}
\Phi(\hat x\cup \hat y)&=&\Phi((\int\hat e\times \hat x)\cup \hat y)\\
&=&\Phi(\int((\hat e\times \hat x)\cup p^*\hat y))\\
&=&\int(\Phi((\hat e\times \hat x)\cup p^*\hat y))\\
&=&\int(\Phi(\hat e\times \hat x)\cup \Phi(p^*\hat y))\\
&=&\int(\Phi(\hat e\times \hat x)\cup p^*\Phi(\hat y))\\
&=&(\int\Phi(\hat e\times \hat x))\cup\Phi(\hat y)\\
&=& \Phi(\int(\hat e\times \hat x))\cup\Phi(\hat y)\\
&=&\Phi( \hat x)\cup\Phi(\hat y)\ . 
\end{eqnarray*}
\hB

\section{The flat theory}\label{uciacascc}

Let $(\hat E,R,I,a)$ be a smooth extension of a generalized cohomology theory $E$.
Let $\hat x\in \hat E^{*+1}([0,1]\times M)$ and $f_0,f_1\colon M\to [0,1]\times M$
be induced by the inclusions of the endpoints. 
Then we have the following homotopy formula.
\begin{lem}\label{udqwdqwdqw}
We have $$f_1^*\hat x-f_0^*\hat x=a\left(\int_{[0,1]\times M/M} R(\hat x)\right)\ .$$
\end{lem}
\proof
Let $p\colon [0,1]\times M\to M$ denote the projection. Then there exists a form
$\rho\in \Omega^{*}([0,1]\times M,\tV)$ such that $\hat x=p^*f_0^*\hat x+a(\rho)$.
We have $R(\hat x)=p^*f_0^*R(\hat x)+d\rho$. Let us write $\rho=\rho_0+dt\wedge \rho_1$ for $t$-dependent forms $\rho_0,\rho_1$  on $M$, where $t$ is the coordinate of the interval $[0,1]$.
We get
$$i_{\partial t} R(\hat x)=\partial_t \rho_0-d\rho_1\ .$$
Integrating we get
$$\int_{[0,1]\times M/M}R(\hat x)=(\rho_0)_{|t=1}-(\rho_0)_{|t=0}-d\int_0^1 \rho_1dt\ .$$
We get
$$f_1^*\hat x-f_0^*\hat x=a(f_1^* \rho-f_0^*\rho)=a((\rho_0)_{|t=1}-(\rho_0)_{|t=0})=a\left(\int_{[0,1]\times M/M}R(\hat x)\right)\ .$$
\hB

\begin{ddd}\label{uiqwdqwdwqd}
We define the \textbf{flat theory} as the subfunctor
$$\hat E_{flat}^*(M):=\ker(R\colon \hat E^*(M)\to \Omega_{cl}^*(M,E))\ .$$
\end{ddd}
Recall that a functor on manifolds is called homotopy invariant if it maps smoothly homotopic smooth maps to the same morphisms. As an immediate consequence of the homotopy formula we get:
\begin{kor}
The flat theory $\hat E_{flat}$ is a homotopy invariant functor.
\end{kor}
As a direct consequence of  \ref{ddd1}.R.3 and \ref{ddd1}.R.1  the flat theory fits into the following natural long exact sequence:
$$\dots\ \stackrel{I}{\to} E^{*-1}(M)\stackrel{\ch}{\to} H^{*-1}(M;\tV)\stackrel{a}{\to} \hat E^*_{flat}(M)\stackrel{I}{\to} E^*(M)\stackrel{\ch}{\to } H^{*}(M;\tV)\stackrel{a}{\to}\dots \ .$$
There is a natural topological candidate for the flat subfunctor which fits into a similar sequence, see (\ref{ufifwefwefewf}).

Recall the construction of the Moore spectrum $\bM G$ for an abelian group $G$. We choose a free resolution
$$0\to \bigoplus_{v\in V} \Z v\stackrel{\alpha}{\to} \bigoplus_{w\in W} \Z w\to G\to 0$$
for  suitable sets $V,W$. Then we define a map of spectra
$$\hat \alpha\colon \bigvee_{v\in V}\bS \to \bigvee_{w\in W} \bS$$ 
which realizes $\alpha$ in reduced integral homology, where we identify
$H\Z_*(\bigvee_{v\in V}\bS)\cong  \bigoplus_{v\in V}\Z v$ and $ H\Z
_*(\bigvee_{w\in W}\bS)\cong  \bigoplus_{w\in W}\Z w$. The isomorphism class of the Moore spectrum
$\bM G$ is defined to fit into the distinguished triangle of the stable homotopy category
$$\bigvee_{v\in V}\bS \stackrel{\hat \alpha}{\to} \bigvee_{w\in W} \bS\to \bM
G\to \Sigma \bigvee_{v\in V}\bS\ .$$
Note that we can and will take $\bM \Z:=\bS$.
We fix an element $1\in G$.
We assume that there is one generator $w_0\in W$ which maps to $1\in G$.
Then we let $\bM \Z\to \bM G$ be the map given by the composition
$\bS\to  \bigvee_{w\in W} \bS\to \bM G$, where the first map is the inclusion of the component
with label $w_0$.

For a spectrum $\bE$ we define
$\bE G:=E\wedge \bM G$. There is a natural identification
$\bE\cong \bE \Z$ and an induced morphism
$\bE\to \bE \R$. The spectrum $\bE\R$ also represents a cohomology theory which admits a canonical isomorphism
$$E\R^*(X)\cong H^*(X;\tV)\ .$$
  
We extend the natural map $\bE\to \bE\R$ to an exact triangle 
$$\Sigma^{-1} \bE \R/\Z\to \bE\to \bE\R\to \bE\R/\Z$$
thus defining a spectrum $\bE\R/\Z$. Note that $\bE\R/\Z\cong \bE\wedge \bM\R/\Z$ so that our notation is consistent.
The fibre sequence  induces a long exact sequence in cohomology
\begin{equation}\label{ufifwefwefewf}
\dots\ \to E^{*-1}(M)\to E \R^{*-1}(M)\to E\R/\Z^{*-1}(M)\to E^*(M)\to E\R^*(M)\to\dots \ .
\end{equation}
In other words, it is very natural to conjecture  that there is a natural transformation
$\Phi_{flat}\colon \hat E_{flat}^*(M)\to E\R/\Z^{*-1}(M)$ so that the following diagram commutes
\begin{eqnarray}\label{udidwqdwqd}
\xymatrix{E^{*-1}(M)\ar@{=}[d]\ar[r]&E \R^{*-1}(M)\ \ar[r]^\alpha&E\R/\Z^{*-1}(M)\ar[r]& E^*(M)\ar@{=}[d]\ar[r]&E\R^*(M)\\E^{*-1}(M)\ar[r]^\ch&H^{*-1}(M;\tV)\ar[u]^\cong\ar[r]^{a}&\hat E_{flat}^*(M)\ar[r]^I\ar[u]^{\Phi_{flat}}&E^*(M)\ar[r]^{ \ch}&H^*(M;\tV)\ar[u]^\cong}
\end{eqnarray}
Such a transformation automatically would be an isomorphism by the Five Lemma.
\begin{ddd}\label{udiqdwqdwqd54545}
We say that the flat theory is \textbf{topological} if such a natural transformation $\Phi_{flat}$ exists.
\end{ddd}

Recall that, given a cohomology theory represented by a spectrum $\bE$, in
\cite[ Definition 4.34]{MR2192936} Hopkins and Singer have constructed a smooth extension $(\hat
E_{HS},R_{HS},I_{HS},a_{HS})$. Moreover they have shown that $\hat
E_{HS,flat}^*$ is topological. We use the notation
$$\Phi^*_{HS,flat}\colon  \hat E^{*}_{HS,flat}(M)\stackrel{\sim}{\to}E\R/\Z^{*-1}(M)$$ 
for the natural isomorphism coming from \cite[ Equation (4.57)]{MR2192936}.
\begin{theorem}\label{zuddwqdqwdqw}
Assume that the cohomology \textcolor{black}{theory} $E$ is rationally even, and
that $E^m$ is countably generated for all $m\in \Z$.
If $(\hat E,R,I,a,\int)$ is a smooth extension of $E$ with integration  which is defined on all smooth manifolds (or alternatively, on all compact manifolds and $E$ is even with $E^{2m}$ finitely generated for all $m\in \Z$),
then the flat functor $\hat E_{flat}$ is topological.
\end{theorem}
\proof
If we could assume that both theories have an integration then we could employ Theorem \ref{main1}.
Unfortunately an integration for the Hopkins-Singer example has not been worked out yet in detail with all the properties we need. Therefore we follow a different path. We start with Theorem \ref{hjsasas} which gives for all $n\in \Z$
an isomorphism
$$\Phi^{2n}\colon \hat E^{2n}\stackrel{\sim}{\to} \hat E_{HS}^{2n}\ .$$
It restricts to an isomorphism
$$\Phi^{2n}\colon \hat E_{flat}^{2n}\stackrel{\sim}{\to} \hat E_{HS,flat}^{2n}$$
of flat theories.
We get the required $\Phi_{flat}^{2n}\colon \hat E^{2n}_{flat}(M)\to E\R/\Z^{2n-1}(M)$ as the composition
$$\Phi_{flat}^{2n}\colon \hat E_{flat}^{2n}(M)\stackrel{\Phi^{2n}}{\to}\hat E^{2n}_{HS,flat}(M)\stackrel{\Phi^{2n}_{flat,HS}}{\to }E\R/\Z^{2n-1}(M)\ .$$

We want to extend this to odd degrees $2n-1$ so that the diagram
$$\xymatrix{\hat E_{flat}^{2n}(S^1\times M)\ar@/^1cm/@{.>}[rr]_{\Phi^{2n}_{flat}}\ar[r]^{\Phi^{2n}}\ar[d]^\int&\hat E_{HS, flat}^{2n} (S^1\times M)\ar[r]^\cong_{\Phi_{HS,flat}^{2n}}&E\R/\Z^{2n-1}(S^1\times M)\ar[d]^\int\\\hat E_{flat}^{2n-1}(M) \ar@{.>}[rr]^{\Phi_{flat}^{2n-1}}&&E\R/\Z^{2n-2}(M)}$$
commutes.
\begin{lem}\label{zdqwudqwdqw}
For all $m\in \nat$ 
$$\int\colon \hat E_{flat}^{m}(S^1\times M)\to \hat E_{flat}^{m-1}(M)$$ is surjective.
\end{lem}
\proof
Let $\hat x\in \hat E_{flat}^{m-1}(M)$. Then we first consider
$y:=q^*\sigma(I(\hat x))\in E^{m}(S^1\times M)$, where $q\colon S^1\times M\to\Sigma M_+$ is the projection and
$\sigma\colon E^{m-1}(M)\to \tilde E^{m}(\Sigma M_+)$ is the suspension isomorphism.
Since $0=\Rham\circ R(\hat x)=\ch\circ I(\hat x)$ we see that $I(\hat x)$ is a torsion element.
Let $0<N\in \nat$ be such that $NI(\hat x)=0$. We choose a lift
$\hat y\in \hat E^{m}(S^1\times M)$ so that $I(\hat y)=y$. Then $NI(\hat y)=0$ and thus
$N \hat y=a(\rho)$ for some $\rho\in \Omega^{m-1}(S^1\times M,\tV)$.
We now replace $\hat y$ by $\hat y-a(N^{-1}\rho)$. Then still
$I(\hat y)=y$, but in addition $R(\hat y)=0$. 

 We have $I(\int \hat y)=\int I(\hat y )=\int y=I(\hat x)$.
Therefore
$\int \hat y-\hat x=a(\theta)$ for some $\theta \in \Omega^{m-2}(M,\tV)$.
Moreover, $d\theta=R(\int (\hat y))-R(\hat x)=0$. Note that
$\int(dt\times \theta)=\theta$ and $d(dt\times \theta)=0$.
If we further replace $\hat y$ by $\hat y-a(dt\times \theta)$, then
$\int \hat y=\hat x$.
\hB 

We now construct $\Phi^{2n-1}_{flat}\colon \hat
E^{2n-1}_{flat}(M)\textcolor{black}{\to E\R/\Z^{2n-2}(M)}$. Let $\hat x\in \hat E^{2n-1}_{flat}(M)$. Then we choose by Lemma \ref{zdqwudqwdqw} an $\hat y\in \hat E^{2n}_{flat}(S^1\times M)$ such that
$\int(\hat y)=\hat x$. Then we define 
$$\Phi^{2n-1}_{flat}(\hat x):=\int(\Phi^{2n}_{flat}(\hat y))\ .$$
We must show that this is well-defined.
Let $\hat y^\prime\in  \hat E^{2n}_{flat}(S^1\times M)$ be such that $\int \hat y^\prime=\hat x$. Then we must show that
$\int(\Phi^{2n}_{flat}(\hat y))=\int(\Phi^{2n}_{flat}(\hat y^\prime))$.
Let $\hat u=\hat y-\hat y^\prime$. It follows from $\int \hat u=0$, that
$I(\hat u)= p^*v$ for some $v\in E^{2n}(M)$. 
Since $\hat u$ is flat we know that $I(\hat u)$ is torsion. Since $p^*\colon E^{2n}(M)\to E^{2n}(S^1\times M)$ is injective we conclude that $v$ is torsion. Therefore we can choose
a lift $\hat v\in \hat E_{flat}^{2n}(M)$. We further find a form
$\theta\in \Omega^{2n-1}(S^1\times M,\tV)$ such that $p^*\hat v+a(\theta)=\hat u$.
If we apply $R$ to this equality, then we get
  $d\theta=0$. Furthermore, from $\int \hat u=0$ we get
 $a(\int \theta)=0$. Therefore $\Rham(\int \theta)=\ch(z)$ for some $z\in E^{2n-2}(M)$.
We choose $w\in E^{2n-1}(S^1\times M)$ such that $\int(w)=z$.
Then
$$\int \Rham(\theta)=\Rham(\int \theta)=\ch(\int w)=\int \ch(w)\ .$$
We now calculate
\begin{eqnarray*}
\Phi^{2n}_{flat}(\hat u)&=&\int(\Phi^{2n}_{HS,flat}(\Phi^{2n}(p^*\hat v+a(\theta))))\\
&=&\int(\Phi^{2n}_{HS,flat}(p^*\Phi^{2n}(\hat v)+a_{HS}(\theta)))\\
&=&\int(p^*\Phi^{2n}_{HS,flat}( \Phi^{2n}(\hat v))+\Phi^{2n}_{HS,flat}(a_{HS}(\theta)))\\
&=&\int(a(\Rham(\theta)))\\
&=&a(\int(\ch(w)))\\
&=&0\ .
\end{eqnarray*} 
This finishes the proof that $\Phi^{2n-1}_{flat}$ is well-defined.

Let us check that 
$$\xymatrix{H^{2n-2}(M;\tV)\ar[r]^a\ar[d]^\rho_\cong&\hat E_{flat}^{2n-1}(M)\ar[r]^I\ar[d]^{\Phi_{flat}^{2n-1}}&E^{2n-1}(M)\ar@{=}[d]\\
E\R^{2n-2}(M)\ar[r]^\alpha &E\R/\Z^{2n-2}(M)\ar[r]^c&E^{2n-1}(M)}$$
commutes.

First we consider
$x\in H^{2n-2}(M;\tV)$.
We must show that
$$\Phi_{flat}^{2n-1}\circ a(x)= \alpha\circ\rho(x)\ .$$
Let $x=\Rham(\omega)$ for some $\omega\in \Omega^{2n-2}_{cl}(M,\tV)$.
Then we take $dt\times \omega\in \Omega^{2n-1}_{cl}(S^1\times M,\tV)$ so that
$\int (a(dt\times \omega))=a(\omega)$. Therefore
\begin{eqnarray*}
\Phi_{flat}^{2n-1}(a(x))&=&\int(\Phi_{HS,flat}^{2n}(\Phi^{2n}(a(dt\times \omega))))\\
&=&\int(\Phi_{HS,flat}^{2n}(a_{HS}(dt\times \omega)))\\
&=&\int(\alpha(\rho(\Rham(dt\times \omega))))\\
&=&\alpha(\int(\rho(\Rham(dt\times \omega))))\\
&=& \alpha(\rho(\int(\Rham(dt\times \omega))))\\
&=&\alpha(\rho(\Rham(\omega)))\\
&=&\alpha(\rho(x))\ .
 \end{eqnarray*}
Next we consider $\hat x\in \hat E_{flat}^{2n-1}(M)$.
We can choose a lift
$\hat y\in  \hat E_{flat}^{2n}(S^1\times M)$ such that
$\int (\hat y)=\hat x$ and $i^*\hat y=0$.
Then we calculate
\begin{eqnarray*}
c(\Phi^{2n-1}_{flat}(\hat x))&=&
c(\int(\Phi^{2n}_{HS,flat}(\Phi^{2n}(\hat y))))\\
&=&\int(c(\Phi^{2n}_{HS,flat}(\Phi^{2n}(\hat y))))\\
&=&\int(I_{HS}(\Phi^{2n}(\hat y)))\\
&=&\int (I(\hat y))\\&=&I(\hat x)\ .
\end{eqnarray*}

In Section \ref{dlede} we will show independently \textcolor{black}{of} the Hopkins-Singer construction that
the flat theory of a smooth extension with integration gives rise to a generalized cohomology theory
which can be compared to $E\R/\Z$.

\hB

\section{Exotic additive structures}

In Theorem \ref{main1} we have seen that a smooth extension of a rationally even cohomology theory
with integration is unique up to unique isomorphism.  In this section we show by example that 
if one disregards the integration there might be many different smooth extensions.

As an example we dicuss $K$-theory. This cohomology theory is even. The
 even part of a smooth extension is unique up to isomorphism because of
Theorem \ref{hjsasas}. In the present section we show that one can change the additive structure of $\hat K^1$ in a non-trivial way.
For simplicity we consider $K$-theory and its smooth extension as two-periodic theories. 
We start with a first example.  Then we show that there are actually infinitely many non-equivalent  additive structures on $\hat K^1$.

Let us start with some smooth extension $(\hat K,R,I,a)$ of complex $K$-theory with additive structure denoted by $+$.
We define a new additive structure on $\hat K^1(M)$ by 
$$\hat u*\hat v:=\hat u+\hat v+a(\frac{1}{2}\ch(I(\hat u)\cup I(\hat v)))\ ,\quad \hat u,\hat v\in \hat K^1(M)\ .$$
Let us verify associativity and commutativity.
\begin{eqnarray*}
(\hat u*\hat v)*\hat w&=&(\hat u+\hat v+a(\frac{1}{2}\ch(I(\hat u)\cup I(\hat v))))*\hat w\\
&=&\hat u+\hat v+a(\frac{1}{2}\ch(I(\hat u)\cup I(\hat v))+\hat w+a(\frac{1}{2}\ch(I(\hat u*\hat v)\cup I(\hat w))))\\
&=&\hat u+\hat v+\hat w+a(\frac{1}{2}\left(\ch(I(\hat u)\cup I(\hat v))+\ch((I(\hat u)+I(\hat v))\cup I(\hat w))\right))\\
&=&\hat u+\hat v+\hat w+a(\frac{1}{2}\ch(I(\hat u)\cup I(\hat v) +I(\hat u)\cup I(\hat w)+I(\hat v)\cup I(\hat w)))\\
&\dots&\\&=&
\hat u*(\hat v*\hat w)\\
\hat u*\hat v&=&\hat u+\hat v+a(\frac{1}{2}\ch(I(\hat u)\cup I(\hat v)))\\
&=&\hat v+\hat u+a(\frac{1}{2}\ch(I(\hat v)\cup I(\hat u)))+a(\ch(I(\hat u)\cup I(\hat v)))\\
&=&\hat v*\hat u\ .
\end{eqnarray*}
Observe that this new additive structure is still compatible with the structure maps
$R,I,a$. In fact, the additional term $a(\frac{1}{2}\ch(I(\hat u)\cup I(\hat v)))$ is annihilated by 
$R$ and $I$, and it vanishes if e.g.  $\hat u=a(\omega)$.
Therefore, 
$\hat K^1$ with this new additive structure together with the old $\hat K^0$ and the old structure maps gives rise to a smooth extension $(\hat K^\prime,R,I,a)$ of $K$-theory.


\begin{prop}\label{jdhbwqdqwdwqdwqdwd}
The smooth extensions  $(\hat K ,R,I,a)$ and  $(\hat K^\prime,R,I,a)$ are not equivalent.
\end{prop}
\proof
Assume that there was  a natural isomorphism
$\Phi\colon \hat K^1\to \hat K^{\prime 1}$ which is compatible with $R$, $I$ and $a$ and such that 
\begin{equation}\label{zduqwdqwd}
\Phi(u+v)=\Phi(u)*\Phi(v)
\end{equation}
Note that $K^1=K^{\prime 1}$ as set-valued functors.
We define $\hat \delta\colon \hat K^1\to \hat K^1$ by
$$\Phi(\hat u)= \hat u+\hat \delta(\hat u)\ .$$
Since $\Phi$ preserves $R$ and $I$ we see that $\hat \delta$ has values
in  $a(H^{ev}(M;\R)/\im(\ch))$. Furthermore, since $\Phi$ is compatible with
$a$ we conclude that $\hat \delta$ comes from a natural transformation
$$\delta\colon K^1(M)\to H^{ev}(M;\R)/\im(\ch)\ , \quad a(\delta(I(\hat u)))=\hat \delta(\hat u)\ .$$
Equation (\ref{zduqwdqwd}) gives
\begin{equation}\label{udiwqdwqdqd}
\delta(u+v)-\delta(u)-\delta( v)=a(\frac{1}{2}\ch( u\cup v))\ .
\end{equation}
We now consider the smooth manifold $T^2=S^1\times S^1$.
Let $e\in K^1(S^1)\cong \Z$ be the generator.
Let $p,q\colon T^2\to S^1$ be the two projections. Then we define $  u:=p^* e$ and $ v:=q^*  e$.
Note that $H^{ev}(T^2;\R)\cong \R^2$ with basis $\{1,\ori\}$, where $\ori\in H^2(T^2;\R)$ is normalized such that $\langle\ori ,[T^2]\rangle=1$. Then the image of the Chern character is the lattice $\ch(K^0(T^2))= \Z\langle1,\ori\rangle\subset \R\langle1,\ori\rangle$. We have
$\ch(u\cup  v)=\ori$.
In particular, $\frac{1}{2}\ch(u\cup  v)\not\in \im(\ch)$ so that
\begin{equation}\label{zuqicqwcqwcwq}
a(\frac{1}{2}\ch(u\cup v))\not=0\ .
\end{equation} Let $*\in T^2$ be a point.
Then $a(\frac{1}{2}\ch(u\cup v))_{|*}=0$.
We define the smooth map $r:=pq\colon T^2\to S^1$ using the group structure of $S^1$. Then we have the identity of $K$-theory classes
$u+v=r^*e$. 
Since $\delta$ is a natural transformation we furthermore get  $\delta(u)=\delta(p^*e)=p^*\delta(e)$. Note that $\delta(e)\in H^0(S^1;\R)/\im(\ch)\cong \R/\Z$.  
It follows that $\delta(u)=c 1$ for some  constant $c\in \R/\Z$. In the same way $\delta(v)=c 1$ and $\delta(u+v)=c 1$.
Then we have $\delta(u+v)-\delta(u)-\delta( v)=-c1$. 
 It follows from (\ref{udiwqdwqdqd})   by considering the restriction to a point that $c=0\in \R/\Z$.
But then $\delta=0$, and this contradicts (\ref{zuqicqwcqwcwq}). 
\hB

\begin{theorem}\label{udqiwdwqdqwd}
There are infinitely many non-isomorphic smooth extensions of complex $K$-theory.
\end{theorem}
\proof
We start with a smooth extension $(\hat K,R,I,a,\int)$ with integration, e.g the multiplicative smooth extension \cite{bunke-20071} with the integration given by Proposition \ref{uzqwiwqdqd}.
 Since $K$-theory is  even, the associated flat theory is topological by Theorem \ref{zuddwqdqwdqw}, i.e.
there is a natural isomorphism
$$\Phi_{flat}\colon \hat K_{flat}^*(M)\stackrel{\sim}{\to} K\R/\Z^{*-1}(M)\ .$$  
In the following, in order to simplify the notation, we will actually identify the flat $K$-theory with
the $\R/\Z$-K-theory and will not write the isomorphism explicitly.

Different smooth extensions will be obtained by modifying the additive structure on $\hat K^1$.
Any other group structure $*\colon \hat K^1(M)\times \hat
K^1(M)\to \hat K^1(M)$ determines and is determined by a natural
transformation $$\hat B\colon \hat K^1(M)\times \hat K^1(M)\to \hat K^1(M)\ ,\quad
\hat x*\hat y=\hat x+\hat y+\hat B(\hat x,\hat y)\ .$$
Compatibility with $R,I$ and $a$ implies that $\hat B$ comes from a transformation
$$B\colon  K^1(M)\times  K^1(M)\to  H^{ev}(M;\R)/\im(\ch)\ , \quad a(B(I(\hat u),I(\hat v)))=\hat B(\hat u,\hat v)\ . $$
Associativity and commutativity of $*$ are equivalent to the conditions  
\begin{equation}\label{u7qiwdqwdqwd} B(u,v+w)+ B(v,w)= B(u,v)+ B(u+v,w)\end{equation}
and  
\begin{equation}\label{u7qiwdqwdqwd1}B(u,v)=B(v,u)\ .\end{equation}

\begin{def}
\textcolor{green}{  Let $\tZ$ be the group of natural transformations of functors
$$B\colon K^1(M)\times  K^1(M)\to H^{ev}(M;\R)/\im(\ch)$$
 which satisfy the
two conditions \eqref{u7qiwdqwdqwd} and \eqref{u7qiwdqwdqwd1}
for all manifolds $M$ and $u,v,w\in K^1(M)$.}
\end{def}
Given $B\in \tZ$,  the new additive structure $*$ is given by 
$$\hat u*\hat v:=\hat u+\hat v+a(B(I(\hat u),I(\hat v)))\ , \quad \hat u,\hat v\in \hat K^1(M)\ .$$

As in the proof of Proposition  \ref{jdhbwqdqwdwqdwqdwd} we will write a natural transformation
$\Phi\colon   \hat K^1\to \hat K^1$ with 
\begin{equation}\label{udiwqdqwdopwqd21}
\Phi(\hat x+\hat y)=\Phi(\hat x)*\Phi(\hat y)
\end{equation} in the form
$\Phi(\hat x)=\hat x+\hat \delta(\hat x)$ for a natural transformation
$$\hat \delta\colon \hat K^1(M)\to \hat K^1(M)\ .$$ Since $\Phi$ must preserve $R,I$ and is compatible with $a$ we again conclude that $\hat \delta$ comes from
$$\delta\colon K^1(M)\to H^{ev}(M;\R)/\im(\ch)\ ,\quad a(\delta(I(\hat u)))=\hat \delta(\hat u)\ .$$
The equation (\ref{udiwqdqwdopwqd21}) translates to
\begin{equation}\label{zuew}B(u,v)=\delta(u+v)-\delta(u)-\delta(v)\ ,\end{equation}
the analog of (\ref{udiwqdwqdqd}).
Note that $\Phi$ is automatically an isomorphism by the Five-Lemma.

\begin{def}
\textcolor{green}{  Let $\tT\subseteq \tZ$ be the group of transformations of the form
  (\ref{zuew}).  Then the set of isomorphism classes of additive extensions of
  $K^1$ is in bijection with the quotient $\tZ/\tT$.}
\end{def}

For $i,j\in \{0,1,2\}$
let $\pr_{ij}\colon \bK_1\times \bK_1\times \bK_1\to \bK_1\times \bK_1$ be the projection
onto the $(i,j)$ component and $s_{01},s_{12}\colon \bK_1\times \bK_1\times \bK_1\to \bK_1\times \bK_1$ the map which adds the first two or the last two factors using the $H$-space structure representing the additive structure of $K^1(X)$.
Let $G$ be a group valued functor on spaces which can be applied to $\bK_1$ and the products with itself. The examples which we have in mind are
$G(X):=H^{ev}(X;\R)/\im(\ch)$ or $G(X):=H^{ev}(X;\R)$.
The cocycle condition for $r\in G(\bK_1\times \bK_1)$ is defined to be 
\begin{equation}\label{coc}
s_{12}^*r+\pr_{12}^*r=\pr_{01}^*r+s_{01}^*r\ .
\end{equation}
Furthermore, let $F\colon \bK_1\times \bK_1\to \bK_1\times \bK_1$ be the flip.
The symmetry condition for $r\in G(\bK_1\times \bK_1)$ is
\begin{equation}\label{sym}F^*r=r\ .
\end{equation}
\begin{lem}
There is a canonical inclusion
$$i\colon \tZ\hookrightarrow H^{ev}(\bK_1\times \bK_1;\R)/\im(\ch)$$
as the subgroup of all elements which satisfy the   cocycle and symmetry conditions (\ref{coc}) and (\ref{sym}).
\end{lem}
\proof
We consider $H^{ev}(M;\R)/\im(\ch)\subseteq K\R/\Z^0(M)$ in the natural way.
Let $(\cK_i)_{i\in \nat}$ be the approximation of $\bK_1$ as in Proposition \ref{duwidwdwd}. 
The family of products $(\cK_i\times \cK_i)_{i\in \nat}$ is then a similar approximation of $\bK_1\times \bK_1$. We let $U\in K^1(\bK_1)$ be the tautological element.
Let $B\in \tZ$.
The family $(B_i)_{i\in \nat}:=(B(\pr_0^*( x_i^* U),\pr_1^*( x_i^*U)))_{i\in \nat}$
is an element  $(B_i)\in \lim ( K\R/\Z^0(\cK_i\times \cK_i))$.
In view of the Milnor sequence
$$0\to \lim\hspace{-2pt}^1  ( K\R/\Z^{-1}(\cK_i\times \cK_i) )\to
K\R/\Z^0(\bK_1\times \bK_1)\to  \lim ( K\R/\Z^0(\cK_i\times \cK_i))\to 0$$
we can choose a preimage
$\tilde B\in K\R/\Z^0(\bK_1\times \bK_1)$.
This preimage is unique up to phantom elements (i.e. elements coming from the
$\lim^1$-term, compare Definition \ref{dqwudqwdqw}), and therefore unique by
Corollary \ref{uoewqeqwe}.
We have exact sequences, natural in the space $X$,
$$\dots\to K^0(X)\to K\R^0(X)\stackrel{\alpha}{\to} K\R/\Z^0(X)
\stackrel{\beta}{\to} K^1(X)\to \dots\ .$$
By construction we know that $B_i \in K\R/\Z^0(\cK_i\times \cK_i)$
is in the image of $\alpha$. Therefore $\beta((B_i))=0$, and this implies that
$\beta(\tilde B)\in K^1_{phantom}(\bK_1\times \bK_1)\subseteq K^1(\bK_1\times
\bK_1)$ (see Definition \ref{dqwudqwdqw}). From Corollary \ref{hdwqdqwdwq} we
get $\beta(\tilde B)=0$ so that $\tilde B=\alpha (\bar B)$ for some
$\bar B\in K\R^0(\bK_1\times \bK_1)$ which is well-defined up to elements coming from
$K^0(\bK_1\times \bK_1)$. If we apply the Chern character, then we get a well-defined element  
$$i(B)=\ch(\bar B)\in H^{ev}(\bK_1\times \bK_1;\R)/\im(\ch)\ .$$
In this way we construct a map
$$\tZ\to H^{ev}(\bK_1\times \bK_1;\R)/\im(\ch)\ .$$
We now discuss injectivity.
Assume that $i(B)=0$. Then
$\tilde B=0$, and this implies that $(B_i)=0$, i.e. $B_i=0$ for all $i\in \nat$.
We now show that this in turn implies $B=0$.

Indeed, let $M$ be a manifold and $u,v\in K^1(M)$.
Then there exists an $i\in \nat$ and maps $f_u,f_v\colon M\to \cK_i$ such that
$u=f_u^*(x_i^*U)$ and $v=f_v^*(x_i^*U)$. But then by naturality of $B$
$B(u,v)=(f_u,f_v)^*B_i=0$.

We now show that every element  $\tB\in H^{ev}(\bK_1\times \bK_1;\R)/\im(\ch)$ which satisfies
the cocycle and symmetry conditions (\ref{coc}) and (\ref{sym}) is in the image of $i$.
Indeed, $\tB$ induces a natural transformation $B$ by
$B(u,v):=(u,v)^*\tB$ for all compact manifolds $M$ and $u,v\in K^1(M)$ considered here as homotopy classes
of maps $M\to \bK_1$. One checks that
$B$ satisfies the cocycle and symmetry conditions, and that $i(B)=\tB$.
\hB

Let $s\colon \bK_1\times \bK_1\to \bK_1$ be the H-space operation, i.e. the map which represents the additive structure on $K^1(X)$.
\begin{lem}
Under $i\colon \tZ \to H^{ev}(\bK_1\times \bK_1;\R)/\im(\ch)$
the subgroup $\tT\subseteq \tZ$ corresponds to the subgroup of classes of the form
$s^*(x)-\pr_0^*x-\pr_1^*x$ for some $x\in H^{ev}(\bK_1;\R)/\im(\ch)$.
\end{lem}
\proof
We consider a natural transformation
$$\delta\colon K^1(M)\to H^{ev}(M;\R)/\im(\ch)\subseteq K\R/\Z^0(M)$$
of functors on manifolds. The family $(\delta_i)_{i\in \nat}$, 
$\delta_i:=\delta(x_i)\in K\R/\Z^0(\cK_i)$ gives rise to an element 
$(\delta_i)\in \lim ( K^0_{\R/\Z}(\cK_i))$.
It again has a unique lift
$\tilde \delta\in K\R/\Z^0(\bK_1)$
which is mapped to a phantom class under
$K\R/\Z^0 (\bK_1)\to K^1(\bK_1)$. From Corollary \ref{hdwqdqwdwq} we conclude
that it belongs to
the subgroup
$\im(K\R^0  (\bK_1) \to K\R/\Z^0 (\bK_1))\subseteq K\R/\Z^0 (\bK_1)$.
If we apply the Chern character we get a well-defined element
$\bar \delta\in H^{ev}(\bK_1;\R)/\im(\ch)$.
Let $B_\delta$ denote the transformation given by (\ref{zuew}).
Then
$i(B_\delta)=s^*(\bar \delta)-\pr_0^*\bar \delta-\pr_1^*\bar \delta$. 
Commutativity of $s$ means $F^*s^*=s^*$. Furthermore, associativity
of $s$ can be expressed as $s_{01}^*s^*=s_{12}^*s^*$.
These two relations together with some obvious identities for the projections imply that the image of $s^*-\pr_0^*-\pr_1^*$ satisfies the cocycle and symmetry conditions and therefore belongs to $i(\tZ)$.
\hB

Since $\bK$ is a ring spectrum   with free $H_*(\bK_*;\Z)$ we get a Hopf ring
$H_*(\bK_*;\Z)$  which has been calculated e.g. in
\cite{MR1870619}. In particular, the underlying graded ring
$$H_*(\bK_1;\Z)=\Lambda[b_1,b_3,\dots]$$
is an exterior algebra with generators $b_{2k-1}\in H_{2k-1}(\bK_1;\Z)$.
It is free over $\Z$ so that we have by the K\"unneth formula
$H_*(\bK_1\times \bK_1;\Z)\cong H_*(\bK_1;\Z)\otimes_\Z H_*(\bK_1;\Z)$.
We let $$c_{2k+1}\in H^{2k+1}(\bK_1;\Z)\cong \Hom_\Ab(H_*(\bK_1;\Z),\Z)^{2k+1}$$
be dual to $b_{2k+1}$ and so that it annihilates all decomposeables.
Then we know that
$H^*(\bK_1;\Z)\cong \Lambda[[c_1,c_3,\dots]]$.  
Furthermore, we have the K\"unneth formula
$$H^*(\bK_1\times \bK_1)\cong  \Lambda[[c_1,c_3,\dots]]\hat \otimes \Lambda[[c_1,c_3,\dots]]\ .$$
We actually do not have to complete here since the cohomology is finitely generated in every single degree.
The sum $s\colon \bK_1\times \bK_1\to \bK_1$ gives a coproduct
$$s^*\colon H^*(\bK_1;\Z)\to H^*(\bK_1;\Z)\hat\otimes  H^*(\bK_1;\Z)\ ,$$
and by construction the generators $c_i$ are primitive, i.e. they satisfy
$$s^*(c_i)=c_i\otimes 1+1\otimes c_i\ .$$

\begin{lem}\label{udiqwedqw}
If $x\in H^{*}(\bK_1;\Z)$ is primitive, then $x\otimes x\in H^{2*}(\bK_1\times \bK_1;\Z)$ satisfies the cocycle  condition (\ref{coc}).
\end{lem}
\proof 
 If we insert $r:=x\otimes x$ into the left hand side of  (\ref{coc}) we get using primitivity of $x$
$$x\otimes x\otimes 1+x\otimes 1\otimes x+1\otimes x\otimes x\ .$$ 
The right-hand side of (\ref{coc}) yields  
$$x\otimes x\otimes 1+x\otimes 1\otimes x+1\otimes x\otimes x\ .$$
\hB 

The following proposition immediately implies Theorem \ref{udqiwdwqdqwd}.
\begin{prop}\label{dioqwdqwd}
The group $\tZ/\tT$ has infinite order.
\end{prop}
\proof
Recall that $H^*(BU;\Z)\cong \Z[[c_2,c_4,\dots]]$, where we index the Chern classes by their degree. We furthermore have an injection
$H^*(BU;\Z)\subset H^*(BU;\Q)\cong \Q[[c_2,c_4,\dots]]$.
On the integral and rational cohomology of a space $X$ we consider the compatible decreasing filtrations $$F^kH^*(X;\Z):=\prod_{i\ge k}H^i(X;\Z)\ ,\quad F^kH^*(X;\Q):=\prod_{i\ge k}H^i(X;\Z)\ .$$

 Let $\ch_{k}(u)\in H^{k}(X;\Q)$ denote the degree $k$-component of the Chern character
$\ch(u)\in H^{*}(X;\Q)$ of $u\in K^*(X)$.
 If $y\in K^0(BU)$ satisfies $\ch(y)\in F^{2k}H^*(BU;\Q)$, then we know that
$\ch_{2k}(y)\in  H^{2k}(BU;\Z)$. In fact, more is known.
 \begin{lem}
For all $x\in H^{2k}(BU;\Z)$ there exists $u\in K^0(BU)$ such that
$\ch(u)\in F^{2k}H^{*}(BU;\Z)$ and $\ch_{2k}(u)=x$.
\end{lem}
\proof
It has been shown by \cite{MR1443417}
that for all $1\le k\le n$ we have $\ch(\rho_{k,n})\in F^{2k}H^*(BU(k);\Q)$
and $\ch_{2k}(\rho_{k,n})=c_{2k}$, where
$\rho_{k,n}=\lambda^k(\rho_n-n)$, $\rho_n$ is the class in $K^0(BU(n))$ of the universal
bundle over $BU(n)$, and $\lambda^k$ is the $k^{\text{th}}$ $\lambda$-operation.
Under the map $BU(n-1)\to BU(n)$ the restriction of
$\rho_n-n$ is  $\rho_{n-1}-(n-1)$.
Therefore the family of classes
$(\rho_n-n)_{n\in \nat}$ defines an element 
$\rho\in \lim ( K^0(BU(n))\cong K^0(BU))$.
We have
$\ch(\lambda^k\rho)\in F^{2k}H^*(BU;\Q)$ and $\ch_{2k}(\lambda^k \rho)=c_{2k}$.
The class $x\in H^{2k}(BU;\Z)$ can be written as $x=p(c_2,c_4,\dots)$, where $p$ is a homogeneous  integral polynomial in the classes $c_{2i}$.
Using that $\ch$ is a ring homomorphism we see that we can choose $u:=p(\lambda^1\rho,\lambda^2\rho,\dots)$.
\hB 

In the following we take the infinite unitary group $U:=\colim (U(n))$ as a model for $\bK_1$.
The universal bundle
$U\to EU\stackrel{\pi}{\to} BU$ gives rise to a transgression
$T_E\colon \tilde E^*(BU)\to \tilde E^{*-1}(U)$ for every generalized cohomology theory $E$.
Let us describe the transgression geometrically. We  assume that all spaces  have base points.
Let $\bE$ be an $\Omega$-spectrum representing $E$. Let $f\colon BU\to \bE_k$
represent a class  $x\in \tilde E^k(BU)$.
Since $\pi^*x=0$
the composition $f\circ \pi$ admits a zero homotopy
$H\colon I\times EU\to \bE_k$. We identify $U$ with the fibre of $\pi$ over the base point of $BU$. 
Then the restriction of $H$ to $I\times U$ closes up and defines a map
$\Sigma U\to \bE_k$, or equivalently, a map $U\to \Omega \bE_k\cong \bE_{k-1}$.
This latter map represents the class $T_E(x)\in \tilde E^{k-1}(U)$. 
The transgression is natural with respect to transformations of cohomology theories $\phi\colon E^*\to F^*$, i.e we have $\phi\circ T_E=T_F\circ \phi$.

For all $k\in \nat$ we choose classes $u_{2k}\in K^0(BU)$ such
that $\ch(u_{2k})\in F^{2k}H^*(BU;\Q)$ and $\ch_{2k}(u_{2k})=c_{2k}$.
We further define $v_{2k-1}:=T_K(u_{2k})\in K^{-1}(U)\cong K^1(U)$. Since
$\ch\colon K\to H\Q$ is a natural transformation of cohomology theories we have
$\ch(v_{2k-1})=T_{H\Q}\ch(u_{2k})$.
The generators of 
$$H^*(U;\Z)=\Lambda[[c_1,c_3,\dots]]$$
are given by $c_{2k-1}=T_{H\Z}(c_{2k})$.
Note that
$T(F^k H^*(BU;\Q))\subseteq F^{k-1}H^*(U;\Q)$.
We conclude that
$\ch(v_{2k-1})\in F^{2k-1}H^*(U;\Q)$ and $\ch_{2k-1}(v_{2k-1})=c_{2k-1}$,
\textcolor{black}{compare also the appendix of \cite{BS}}.

It is a well-known fact that the transgression $T(x)\in H^{*-1}(U;\Q)$ of  $x\in H^*(BU;\Q)$  is primitive.
In order to see this one could consider a dual transgression in homology.
One easily observes geometrically that the transgression of a non-trivial Pontrjagin
product of homology classes gets annihilated. This implies that $T(x)$ annihilates all non-trivial Pontrjagin products of homology classes. This property is equivalent to primitivity.

We  consider the differential\footnote{The $H$-space $\bK_1$ gives rise to the simplicial space $B^\bullet\bK_1$
and the corresponding cobar complex
$(H^*(B^\bullet\bK_1;\R),d)$ in cohomology. It is for this reason that we call $d$ a differential.}
$d\colon H^*(\bK_1;\R)\to H^*(\bK_1\times \bK_1;\R)$ given by
$$d(x):=s^* x-1\times x-x\times 1=s^*x-\pr_0^*x-\pr_1^*x \ ,$$
where $s\colon \bK_1\times \bK_1
\to \bK_1$ is the $H$-space structure as above.
Note that $H^*(\bK_1;\R)\cong\Lambda_\R[[c_1,c_3,\dots]]$
is generated by primitive elements.
We consider $H^*(\bK_1;\Z)\cong \Lambda_\Z[[c_1,c_3,\dots]]\subset H^*(\bK_1;\R)$
as a subspace in the natural way.
\begin{lem}\label{duqwidqwd}
For $i>0$ the differential $d\colon H^{2i}(\bK_1;\R)\to H^{2i}(\bK_1\times \bK_1;\R)$ is injective.
If $x\in H^{2i}(\bK_1;\R)$ and
$d(x)\in H^{2i}(\bK_1\times \bK_1;\Z)$, then
$x\in H^{2i}(\bK_1;\Z)$.
\end{lem}
\proof
Let $I:=\{i_1<i_2<\dots<i_{r}\}$ be a sequence of odd integers. It determines a monomial
$m_I:=c_{i_1}\dots c_{i_{r}}\in  H^*(\bK_1;\R)$.
Using the primitivity of the $c_i$ we observe that
$$s^* (m_I)=\sum_{P\sqcup Q} l(P,Q) m_P\otimes m_Q\ ,$$
where the sum is taken over all partitions $P\sqcup Q=\{i_1,\dots i_r\}$ and $l(P,Q)\in \{1,-1\}$
is the sign of the permutation determined by $P,Q$. 
It follows that
$$d m_I=\sum_{P\sqcup Q=I, P,Q\not=\emptyset} l(P,Q) m_P\otimes m_Q\ .$$
If $x=\sum_{I\not=\emptyset} a_Im_I\in H^{2i}(\bK_1;\R)$ with $a_I\in \R$, then we have
\begin{equation}\label{iudoqwdqwd}
d(x)=\sum_{I\not=\emptyset}a_I\sum_{P\sqcup Q=I, P,Q\not=\emptyset} l(P,Q) m_P\otimes m_Q\ .
\end{equation}
As $x$ is an even cohomology class, each $I$ in this sum with $a_I\not=0$ has at least two elements so that 
there exists at least one partition
$P\sqcup Q=I$ with $P$ and $Q$ nonempty. We now see that we can recover the $a_I$ from the 
right-hand side of (\ref{iudoqwdqwd}).
In particular, if the right-hand side is integral, all the coefficients $a_I$ must be integral.
\hB

\begin{lem}
The images $w_{4k-2} \in H^{ev}(\bK_1\times \bK_1;\R)/\im(\ch)$ of the classes $$\frac{1}{2} \ch(v_{2k-1})\times \ch(v_{2k-1})\in H^{ev}(\bK_1\times \bK_1;\R)$$ 
generate an infinite sum of copies of $\Z/2\Z$
in $\tZ/\tT$.
\end{lem}
\proof
Since
$\ch(v_{2k-1})=T(\ch(u_{2k}))$ is primitive, the class $w_{4k-2}$ satisfies the cocycle condition
(\ref{coc}) by Lemma \ref{udiqwedqw}. Since $\ch(v_{2k-1})$ is odd and
$\ch(v_{2k-1})\times \ch(v_{2k-1})=\ch(v_{2k-1}\times v_{2k-1})\in \im(\ch)$ and in addition $2w_{4k-2}=0$ we see that $w_{4k-2}$ satisfies the symmetry condition (\ref{sym}). 
It therefore remains to show that the $w_{4k-2}$ are independent.

Assume that
$\sum_{i=1}^r a_i w_{4i-2}=0$. We must show that all $a_i$ are even.
We proceed by induction.
Let us assume that we have shown that $a_i=2b_i$ for $i=1,\dots k-1$.
Then 
$0=\sum_{i=k}^r a_i w_{4i-2}$ so that there are
$x\in K^0(\bK_1\times \bK_1)$ and $y\in H^*(\bK_1;\R)$ such that
$$\sum_{i=k}^r \frac{a_i}{2} \ch(v_{2i-1})\times \ch(v_{2i-1})=d(y)+\ch(x)\ .$$
Note that the left-hand side belongs to
$F^{4k-2}H^*(\bK_1\times \bK_1;\R)$, and that the lowest term is
$\frac{a_k}{2}c_{2k-1}\times c_{2k-1}$. We claim that we can adjust $y$ such that
$d(y)\in F^{4k-2}H^*(\bK_1\times \bK_1;\R)$.
Let $j\in \nat$ be minimal such that
$d(y)\in F^jH^*(\bK_1\times \bK_1;\R)$. 
By the first assertion of Lemma \ref{duqwidqwd} the lowest term of $d(y)$ is given by 
$d(y)_j=d(y_j)$, where $y_j$ is the component of $y$ in degree $j$. If $j<4k-2$, then 
 $\ch(x)\in F^jH^*(\bK_1\times \bK_1;\R)$, and $d(y_j)=-\ch_{j}(x)$. We conclude
that $d(y_j)$ is integral. By the second assertion of Lemma \ref{duqwidqwd}
we see that $y_j$ is integral.
Hence there exists $u\in K^0(\bK_1)$ such that $\ch_j(u)=y_j$.
If we replace $y$ by $y^\prime:=y-\ch(u)$ and $x$ by $x^\prime:=x+d(u)$,
then we have $d(y)+\ch(x)=d(y^\prime)+\ch(x^\prime)$, and we have $d(y^\prime)\in F^{j+1}H^*(\bK_1\times \bK_1;\R)$.
If $j=4k-2$, then
$$\frac{a_k}{2}c_{2k-1}\times c_{2k-1}=d(y_{4k-2})+\ch_{4k-2}(x)\ .$$
Since now $\ch(x)\in F^{4k-2}H^*(\bK_1\times \bK_1;\Q)$
we see that $\ch_{4k-2}(x)$ is integral. In particular, the coefficent of
$\ch_{4k-2}(x)$ in front of the monomial $c_{2k-1}\times c_{2k-1}$ must be integral.
Since this monomial does not occur in $d(y)$ by the calculation (\ref{iudoqwdqwd})
we see that $\frac{a_k}{2}$ must be integral, too.
\hB 
This finishes also the proof of Proposition \ref{dioqwdqwd} and therefore of Theorem \ref{udqiwdwqdqwd}.
\hB

\section{Mayer-Vietoris sequence}\label{dlede}\label{e89wfoewfwefq}

We consider a smooth extension with a natural integration $(\hat
E,R,I,a,\int)$ of a generalized cohomology theory $E$.  It gives rise to the
flat theory $\hat E_{flat}$ (Definition \ref{uiqwdqwdwqd}) which is a homotopy
invariant functor on the category of  manifolds.  In this section we show that
it is a generalized cohomology theory by constructing a Mayer-Vietoris
sequence. If  $\hat E_{flat}$ is topological, then it is clearly a generalized cohomology
theory. The point of the present section is to give a construction of the
Mayer-Vietoris sequence independently of Theorem \ref{zuddwqdqwdqw} and hence
of the Hopkins-Singer theory.

In the present section we assume that the smooth extension is defined on the category of all smooth manifolds.
There is a version of the theory for smooth extensions defined on compact manifolds (with boundary).
In this case we must replace the words manifold by compact manifold and finite-dimensional countable $CW$-complex by finite $CW$-complex at the appropriate places. We indicate some further modifications as footnotes.

We consider a manifold
$M$ which is decomposed as a union of open submanifolds \begin{equation}\label{udiqwdqwdqwdqdqwd}
M=U\cup V\ ,\quad A:=U\cap V\ .
\footnote{The modification in the case of a smooth extension defined on compact manifolds is as follows. We assume that $U$ and $V$ are closed (with boundary), and that there are deformation retracts  of $U$ and $V$ onto compact $U^\prime\subset \inter(U)$ and $V^\prime\subset \inter V$.}\end{equation}

We choose a smooth function
$\chi\colon M\to [-1,1]$ such that
$$\chi_{M\setminus V}=-1\ ,\quad \chi_{|M\setminus U}=1\footnote{If we must work with compact manifolds, then we require these condition with $U^\prime,V^\prime$ in place of $U,V$.}\ .$$
Let $\Sigma^uA:=[-1,1]\times A/\sim$ denote the unreduced suspension, where
$(t,a)\sim (t^\prime,a^\prime)$ if and only if $t=t^\prime=1$ or
$t=t^\prime=-1$ or $t=t^\prime$ and $a=a^\prime$.
The equivalence classes in the first two cases will be denoted by $*_+$ and $*_-$.
We define a projection
$p\colon M\to \Sigma^u A$ by 
\begin{equation}\label{zuiudwqdqwdqd64646}
p(m)=\left\{ \begin{array}{ccc}
(\chi(m),m)&m\in A\\
*_-&m\in U\setminus A\\
*_+&m\in V\setminus A
\end{array}\right\}\footnote{If we must work with compact manifolds, then the modified formula would involve the retractions of $U$ \textcolor{black}{onto} $U^\prime$ and $V$ \textcolor{black}{onto} $V^\prime$.}
\end{equation}
The restriction  $p_{|U}$ is zero-homotopic.
Let us give the zero homotopy
$p_U\colon I\times U\to \Sigma^uA$ by an explicit formula:
$$p_U(t,m):=(t+(1-t)\chi(m),m)\ .$$
We define the zero homotopy $p_V\colon I\times V\to \Sigma^u A$ of $p_{|V}$ by a similar formula.
The decomposition (\ref{udiqwdqwdqwdqdqwd}) gives rise to a Mayer-Vietoris sequence
$$\dots \to E^k(M)\to E^k(U)\oplus E^k(V)\to E^k(A)\stackrel{\partial}{\to} E^{k+1}(M)\to \dots\ .$$ 
An explicit description of the boundary operator is as follows.
We consider a class $x\in E^k(A)$ which we assume to be represented by a map
$g\colon A\to \bE_k$.
Then using the projection $\Sigma^u\bE_k\to \Sigma \bE_k$ and the structure map
$\Sigma\bE_k\to \bE_{k+1}$ we define
$\Sigma g\colon \Sigma^u A\to \Sigma^u \bE_k\to \Sigma \bE_k \to \bE_{k+1}$.
The composition
$h:=\Sigma g\circ p\colon M\to \bE_{k+1}$ represents
$\partial x\in E^{k+1}(M)$.
It comes with zero homotopies
$h_U:=\Sigma g\circ p_U$ and $h_V:=\Sigma g\circ p_V$ of $h_{|U}$ and $h_{|V}$.
Over $A$ we can glue these zero homotopies
to a map
$$h_U^{op}\sharp h_V\colon S^1\times A\cong \frac{[-1,0]\times A\sqcup [0,1]\times A}{\sim}\stackrel{ h_U^{op}\sqcup h_V}{\to} \bE_{k+1}\ ,$$
where $h_U^{op}(t,m)=h_U(-t,m)$. 
This map represents a class
$\tilde x\in E^{k+1}(S^1\times A)$ such that $\int\tilde x=x$.


\begin{lem}\label{dwqdqwdw}

Let $x\in E^k(M)$ be such that $x_{|U}=0$ and $x_{|V}=0$. Then there exists
a based manifold $N$,   a class $y\in E^k(N)$ with trivial restriction to the base point, a map $f\colon M\to N$ such that
$x=f^*y$,  and zero homotopies $f_U$ and $f_V$ of 
$f_{|U}$ and $f_{|V}$.

In addition,
if $z\in E^{k-1}(A)$ is such that $x=\partial z$ then we can choose
$f\colon M\to N$, $y\in E^k(N)$  and the zero homotopies  
$f_U$ and $f_V$ such that $\int l^*y=z$, where
$l:=f_U^{op}\sharp f_V$.

Even more specifically,
if $\partial z=x=0$, then we can choose these objects such that in addition there exists a zero homotopy $f_M$ of $f$.
\end{lem}
\proof
Let $g\colon M\to \bE_k$ represent the class $x$. Then there are zero homotopies
$g_U\colon I\times U\to \bE_k$ and $g_V\colon I\times V\to \bE_k$ of $g_{|U}$ and $g_{|V}$.
Since $M$, $I\times U$ and $I\times V$ are finite-dimensional manifolds there
exists a countable finite-dimensional subcomplex
$X\subseteq \bE_k$ containing the images of the maps $g,g_U$ and $g_V$.
Every countable finite-dimensional $CW$-complex is homotopy equivalent to a smooth manifold.
We choose such a manifold $N$ with base point together with  mutually inverse homotopy equivalences
$s\colon N\to X$ and $r\colon X\to N$,  and we define $y\in E^k(N)$ as the class represented by the composition
$N\stackrel{s}{\to} X\to \bE_k$. We further set $\tilde f:=r\circ g$, 
$\tilde f_U:= r\circ g_U$, and $\tilde f_V:= r\circ g_V$.
Finally we first replace $\tilde f$ by a homotopic smooth map. In the second step we replace
the zero homotopies $\tilde f_U, \tilde f_V$ by smooth zero homotopies  $f_U,f_V$ of $f_{|U} , f_{|V}$. 

If $x=\partial z$, then we can choose a map 
$j\colon A\to \bE_{k-1}$ representing $z$ and the map $g$   as the composition  $$g\colon M\stackrel{p}{\to} \Sigma^u A\stackrel{\Sigma j}{\to} \Sigma^u \bE_{k-1}\to \bE_k\ .$$
In this case the homotopies $g_U,g_V$ are induced by the homotopies $p_U,p_V$.
Then we proceed as above.

Finally, if $x=0$, then the map $g$ admits a zero homotopy $g_M$.
We choose the finite subcomplex $X$ sufficiently large to capture the image of $g_M$, too.
Then we proceed as above and let $f_M$ be induced by $g_M$.
\hB 

\begin{lem}\label{duiqdwqdqwd}
Let $u\in E^{k-1}(U)$. Then there exists a based manifold $N$, a class $y \in E^{k}(N )$ vanishing on the base point, and a  map $l \colon \Sigma^u U\to N $ such that $l^*y=\sigma(u)$. 
We can assume that $l $ is constant near the singular points of the unreduced suspension and smooth elsewhere, \textcolor{black}{and} that the singular points $*_\pm$ are mapped to the base point $*$ of $N $.
\end{lem}
\proof Let $\bE_k$ be as in the proof of  Lemma \ref{dwqdqwdw}. 
We choose a map $g\colon \Sigma^u U\to \bE_k$ which represents the class $\sigma(u)$ and maps $*_\pm$ to the base point.
The map $g$ factors over a finite-dimensional countable subcomplex $X$ which we approximate by a smooth manifold $s\colon N\to X$, $r\colon X\to N$ such that the composition $N\stackrel{s}{\to} X\to \bE_k$ represents $y\in E^k(N)$
 and maps the basepoint $*\in N$ to the base point of $\bE_k$.
We set $\tilde l :=r¸\circ g$. Finally we replace this by a homotopic map with the required properties.
\hB 

\begin{lem}\label{7eifdwqedqfqd}
Let $u\in E^{k-1}(U)$ and $v\in E^{k-1}(V)$ be such that $u_{|A}=v_{|A}$.
Then there exists a based smooth manifold $N$, a class $y\in E^k(N)$ vanishing on the base point, and maps $f\colon U\to N$ and $g\colon V\to N$ such that $u=f^*y$, $v=g^*y$, and there is a homotopy $f_{|A}\sim g_{|A}$.
\end{lem}
\proof
 Let $\bE_k$ be as in the proof of  Lemma \ref{dwqdqwdw}. Furthermore, let
 $a\colon U\to \bE_k$ and $b\colon V\to \bE_k$ represent
the classes $u$ and $v$.
Then there exists a homotopy $h\colon a_{|A}\sim b_{|A}$.
We choose a finite-dimensional countable $X\subseteq\bE_k$ over which the maps
$a,b$ and the homotopy $h$ factor. Then we choose a smooth approximation
$s\colon N\to X$, $r\colon X\to N$ by homotopy equivalences, let $y\in
E^k(N)$ be represented by $N\stackrel{s}{\to} X\to \bE_k$,
$\tilde f:=r\circ a$, $\tilde g:=r\circ b$, and $\tilde h:=r\circ h$.
Then we first approximate $f$ and $g$ by smooth maps, and then choose the smooth
homotopy between $f_{|A}$ and $g_{|A}$ by adapting $\tilde h$.
\hB

Let us now construct the boundary operator
$$\hat \partial\colon \hat E_{flat}^k(A)\to \hat E^{k+1}_{flat}(M)\ .$$
Let $\hat x\in \hat E_{flat}^k(A)$ be given. We set $x:=I(\hat x)$. Then by Lemma \ref{dwqdqwdw} there is
a pointed  smooth manifold $N$, a class $y\in E^{k+1}(N)$ with $y_{|*}=0$, and a smooth map
$f\colon M\to N$ such that
\begin{enumerate}
\item $f^*y=\partial x$,
\item there are zero homotopies $f_{U}\colon I\times U\to N$ and $f_V\colon I\times V\to N$
of $f_{|U}$ and $f_{|V}$,
\item with  $l:=f_U^{op}\sharp f_V\colon S^1\times A\to N$ we have $\int l^*y=x$.
\end{enumerate}
We choose a smooth lift $\hat y\in \hat E^{k+1}(N)$ which restricts
trivially to the base point. Indeed, if we choose some class
$\hat y\in \hat E^{k+1}(N)$ such that $I(\hat y)=y$, then       
$\hat y_{|*}=a(c)$ for some $c\in \Omega^k(*,\tV)\cong  \tV^k $. We denote the constant zero form with value
$c$ on $N$ by the same symbol. If we replace $\hat y$ by $\hat y-a(c)$, then the  restriction of the  new $\hat y$ to the base point vanishes.

We furthermore choose
  a form
$\rho\in \Omega^{k-1}(A,\tV)$ such that 
$$\int l^*\hat y+a(\rho)=\hat x\ .$$
Finally, we   choose forms $\rho_U\in  \Omega^{k-1}(U,\tV)$ and
$\rho_V\in  \Omega^{k-1}(V,\tV)$ such that
$\rho=\rho_{V|A}-\rho_{U|A}$.
Using these choices we define a form
$\kappa\in \Omega^{k}(M,\tV)$
by the description
$$\kappa_{|U}:=\int_{I\times U/U} f_U^* R(\hat y)+d\rho_U\ ,\quad \kappa_{|V}:=\int_{I\times V/V} f_V^* R(\hat y)+d\rho_V\ .$$
Note that $\kappa$ is well-defined on $A$ since
\begin{eqnarray*}\lefteqn{
 \left(\int_{I\times V/V} f_V^* R(\hat y)+d\rho_V\right)_{|A}- \left(\int_{I\times U/U} f_U^* R(\hat y)+d\rho_U\right)_{|A}}&&\\&=&
\int_{S^1\times A/A} l^* R(\hat y)+d\rho\\&=&R(\hat x)\\&=&0\ .
\end{eqnarray*}
We define
$$\hat \partial \hat x:=f^*\hat y+a(\kappa)\ .$$

\begin{lem}
$\hat \partial \hat  x$ is well-defined.
\end{lem}
\proof
We show step by step that if we alter the choices going into the construction of $\hat \partial\hat x$
we get the same result. We indicate the changed objects by a prime.
\begin{enumerate}
\item
If we choose $\rho_U^\prime$ and $\rho_V^\prime$ such that
$\rho^\prime_{V|A}-\rho^\prime_{U|A}=\rho$, then
there exists a form $\theta\in  \Omega^{k-1}(M,\tV)$ such that
$\theta_{|U}=\rho_U-\rho_{U}^\prime$ and $\theta_{|V}=\rho_V-\rho_{V}^\prime$.
We thus get $\kappa^\prime=\kappa-d\theta$ and hence
$$\hat \partial^\prime \hat x=f^*\hat y +a(\kappa^\prime)=f^*\hat y +a(\kappa)-a(d\theta)=f^*\hat y +a(\kappa)=\hat \partial \hat x\ .$$
\item Let us choose another form $\rho^\prime$ such that
$\int l^*\hat y+a(\rho^\prime)=\hat x$.
Then $\rho^\prime=\rho+\theta$ with $d\theta=0$ and
$\Rham(\theta)=\ch(u)$ for some $u\in E^{k-1}(A)$.
As for $\rho$ we choose a decomposition $$\theta_U\in \Omega^{k-1}(U,\tV)\ ,\quad  \theta_V\in \Omega^{k-1}(V,\tV)\ ,\quad  \theta_{V|A}-\theta_{U|A}=\theta\ .$$   Furthermore,
we set $\rho_U^\prime:=\rho_U+\theta_U$ and $\rho_V^\prime:=\rho_V+\theta_V$.
Then $\kappa^\prime=\kappa+\lambda$, where $\lambda\in \Omega^k(M,\tV)$ is the
closed form determined by $\lambda_{|U}=d\theta_U$ and
$\lambda_{|V}=d\theta_V$. Its cohomology  class is given by
$$\Rham(\lambda)=\partial \Rham(\theta)=\partial \ch(u)=\ch(\partial  u)\ .$$ Therefore  we have
$a(\kappa)=a(\kappa^\prime)$.
This implies $\hat \partial^\prime \hat x=\hat \partial \hat  x$.
\item
If we choose another  lift $\hat y^\prime$, then
$\hat y^\prime=\hat y+a(\theta)$ for
$\theta\in  \Omega^{k}(N,\tV)$.
We get
$$\int_{I\times U/U} f_U^*R(\hat y^\prime)=\int_{I\times U/U}
\left(f_U^*R(\hat y)+ f_U^*d\theta\right)=\int_{I\times U/U} f_U^*R(\hat y)
-d\textcolor{black}{(}\int_{I\times U/U}f_U^*\theta\textcolor{black}{)} -f_{|U}^*\theta\ .$$
A similar formula holds true over $V$.
Note that $\int l^*\hat y^\prime-\int l^*\hat y=a(\int l^*\theta)$ and therefore
$\rho^\prime-\rho= \int l^*\theta$.
We can take
$\rho_{U}^\prime:=\rho_U+\int_{I\times U/U} f_U^*\theta$ and a similar formula over $V$. 
We see that
$\kappa^\prime=\kappa- f^*\theta$ and hence
$$\hat \partial^\prime\hat x=f^*\hat y^\prime +a(\kappa^\prime)=f^*\hat y+a(f^*\theta) +a(\kappa-f^*\theta)=f^*\hat y+a(\kappa)=\hat \partial \hat x\ .$$
\item\label{ttwer2}
Let $N^\prime$ be a   smooth manifold with class $y^\prime\in  E^{k+1}(N)$ and smooth map $u\colon N\to N^\prime$ such that $u^* y^\prime= y$.
If $f^\prime=u\circ f$ and $f_U^\prime=u\circ f_U$ and $f_V^\prime=u\circ f_V$, and we choose
a lift $\hat y^\prime\in \hat E^{k+1}(N^\prime)$ and take $\hat y:=u^*\hat y^\prime$,
then we   get
$ \hat \partial^\prime\hat x=\hat \partial \hat x$.
\item\label{ttwer1}
We now assume that we have a fixed manifold $N$ with class $y\in E^{k+1}(N)$
and homotopic choices $f\sim f^\prime$ and compatible homotopic choices of homotopies
$f_U\sim f_U^\prime$ and $f_V\sim f_V^\prime$.

We then consider  the decomposition  $I\times M=I\times U\cup I\times V$.
The homotopies give maps
$F\colon I\times M\to N$ and zero homotopies $F_U\colon I\times I\times U\to N$ and $F_V\colon I\times I\times V\to N$.  
  Applying the construction of the boundary operator in this case gives a class
$\hat \partial (\pr_A^*\hat x)\in \hat E^{k+1}_{flat}(I\times M)$ which restricts to
the classes
$\hat \partial\hat x$ and $\hat \partial^\prime\hat x$  at the two boundary components. By the homotopy invariance  of the functor $E^{k+1}_{flat}$ we conclude that
 these two classes coincide.
\item Finally assume that  we have two choices
$(N_j,y_j,f_j,f_{j,U},f_{j,V})$, $j=0,1$. Then by a similar argument
as in the proof of Lemma \ref{dwqdqwdw} there exists a third
choice
$(N,y,f,f_{U},f_{V})$ together with maps
$u_j\colon N_j\to N$ and homotopies
$u_j^*f\sim f_j$ and compatible homotopic choices  of homotopies $u^*_jf_{U}\sim f_{j,U}$ and $u^*_jf_{V}\sim f_{j,V}$.
By a combination of \ref{ttwer1} and \ref{ttwer2} we conclude that 
the result for $\hat \partial \hat x$ constructed with $j=0$ coincides with the result for $j=1$.
\end{enumerate}
 \hB

\begin{lem}
The boundary operator
$\hat \partial\colon E^k_{flat}(A)\to E^{k+1}_{flat}(M)$
is additive.
\end{lem}
\proof
Let $\hat x_i\in E^k_{flat}(A)$, $i=0,1$ be given.
We do the construction of $\hat \partial \hat x_i$ based on the choice
$f_i\colon M\to N_i$, $\hat y_i\in \hat E^{k+1}(N_i)$.
Then we can use the choice
$N:=N_0\times N_1$ with the class
$\hat y:=\pr_0^*\hat y_0+\pr_1^*\hat y_1$, the map
$f:=(f_0,f_1)$ and the zero homotopies
$f_{U}=(f_{0,U},f_{1,U})$ and $f_{V}=(f_{0,V},f_{1,V})$ for the construction of $\hat \partial (\hat x_1+\hat x_1)$.
If we have choosen
$\rho_{i,U},\rho_{i,V}$, then $\rho_U:=\rho_{0,U}+\rho_{1,U}$
and $\rho_V:=\rho_{0,V}+\rho_{1,V}$ is an appropriate choice for the sum.
We get $\kappa=\kappa_0+\kappa_1$ and finally
$\hat \partial \hat x=\hat \partial \hat x_0+\hat \partial \hat x_1$.
 \hB

\begin{lem}\label{zqidqwdwqdqw}
The boundary operator is natural.
\end{lem}
\proof
Let $M^\sharp =U^\sharp\cup V^\sharp$ be a decomposition of another manifold with a smooth map
$\phi\colon M^\sharp\to M$ such that $\phi(U^\sharp)\subseteq U$ and $\phi(V^\sharp)\subseteq V$.
Let $A^\sharp:=U^\sharp\cap V^\sharp$.
If $\hat x\in \hat E_{flat}^k(A)$, then we must show that
$$\hat \partial^\sharp \phi_{|A^\sharp}^*\hat x=\phi^*\hat \partial \hat x\ .$$
We use the notation introduced in the construction of $\hat \partial$.

We choose $f^\sharp:=f\circ \phi\colon M^\sharp\to N$.
Then $$f^\sharp y=\phi^*f^*y =\phi^*\partial x=\partial \phi^*_{|A^\sharp} x\ .$$
For the zero homotopies we choose $f^\sharp_U:=f_U\circ \phi_{|U^\sharp}$ and
$f^\sharp_V:=f_V\circ \phi_{|V^\sharp}$. Then for the loop we get
$l^\sharp=l\circ (\id_{S^1}\times \phi_{|A^\sharp})\colon S^1\times A\to N$.
It follows that
$$\int l^{\sharp,*} y=\int (\id_{S^1}\times \phi_{|A^\sharp})^* l^*y=
\phi_{|A^\sharp}^* \int l^*y=\phi_{|A^\sharp}^* (x)\ .$$
Hence we construct $\hat \partial^\sharp$ using $f^\sharp\colon M^\sharp\to N$, the homotopies $f^\sharp_U,f^\sharp_V$ and the class
$\hat y$. 
Indeed, we can choose $\rho^\sharp:=\phi_{|A^\sharp}^*\rho$ since
$$\int l^{\sharp,*} \hat y+a(\phi_{|A^\sharp}^*\rho)=
\phi_{|A^\sharp}^*\int l^* \hat y +\phi_{|A^\sharp}^* a(\rho)=\phi_{|A^\sharp}^*\hat x\ .$$
Furthermore we can take
$\rho^\sharp_U:=\phi_{|U^\sharp}^*\rho_U$ and $\rho^\sharp_V:=\phi_{|V^\sharp}^*\rho_V$
which leads to $\kappa^\sharp=\phi^*\kappa$. After all this we see that
$$\hat \partial^\sharp \phi_{|A^\sharp}^*\hat x=\phi^*\hat \partial \hat x\ .$$ \hB

\begin{prop}
We have a long exact Mayer-Vietoris sequence
$$\dots \hat E_{flat}^k(M)\stackrel{c}{\to} \hat E_{flat}^k(U)\oplus \hat E_{flat}^k(V)\stackrel{b}{\to} \hat E_{flat}^k(A)\stackrel{\hat \partial}{\to} \hat E_{flat}^{k+1}(M)\to \dots\ .$$ 
\end{prop}
\proof
The map $c$ is given by
$$c(\hat x):=\hat x_{|U}\oplus \hat x_{|V}\ .$$ The map $b$ is defined by
$$b(\hat x\oplus\hat y):=\hat y_{|A}-\hat x_{|A} \ .$$
We first show that this is a complex.
\begin{enumerate}
\item
It is clear that $b\circ c=0$.
\item \label{uiwdqwdqw}
Next we show that $c\circ \hat \partial=0$.
We have, using the homotopy formula Lemma \ref{udqwdqwdqw} and the vanishing of $\hat y$ on the base point,
$$\hat \partial \hat x_{|U}=(f^*\hat y)_{|U} +a(\kappa)_{|U}=-a(\int_{I\times U/U} f_U^*R(\hat y))+a(\kappa_{|U})=a(d\rho_U)=0\ .$$
\item
Next we verify that $\hat \partial \circ b=0$.
Let for example $\hat x\in \hat E^{k}_{flat}(U)$ be given. Then we set
$U^\sharp:=U$ and $V^\sharp=M$ and consider the second decomposition
$M=U^\sharp\cup V^\sharp$ with $A^\sharp=U$.
The map $\phi=\id\colon M\to M$ respects these decompositions.
By \ref{uiwdqwdqw}. we have \textcolor{black}{on the one hand}
$\hat \partial^\sharp \hat x=(\hat \partial^\sharp \hat x)_{|V^\sharp}=0$.
On the other hand, by Lemma \ref{zqidqwdwqdqw}
$$\hat \partial\hat x=\hat \partial\phi^*\hat x=\phi^*\hat \partial^\sharp \hat x=0\ .$$
 \end{enumerate}
 We now verify exactness.
\begin{enumerate}
\item Let $\hat x\in \hat E^{k}_{flat}(A)$ be such that $\hat \partial \hat x=0$.
In this case $f^*\hat y+a(\kappa)=0$. By Lemma \ref{dwqdqwdw} we can in addition assume that there is a zero homotopy $f_M$. We get by the homotopy formula Lemma \ref{udqwdqwdqw} and the vanishing of $\hat y$ on the base point
$$0=-a(\int_{I\times M/M} f_M^* R(\hat y))+a(\kappa)\ .$$
Hence there is a class $\hat m\in \hat E^{k-1}(M)$ such that
\begin{equation}\label{uoiqdqwdqwwqd}
\kappa-\int_{I\times M/M} f_M^* R(\hat y)=R(\hat m)\ .
\end{equation}
First we will make a modification which allows us to assume that $R(\hat m)=0$.
To this end by Lemma \ref{duiqdwqdqwd} we choose a manifold $N_M$, a class $y_M\in E^{k-1}(N_M)$, and a  smooth map $l_M\colon \Sigma^u M\to N_M$ such that $l_M^*y_M=-\sigma(I(\hat m))$. 
We can assume that $l_M$ is constant near the singular points and smooth elsewhere.
We choose a smooth lift $\hat y_M\in \hat E^{k-1}(N_M)$. Then we set
$N^\prime:=N\times N_M$ and $\hat y^\prime:=\pr_N^*\hat y+\pr_{N_M}^*\hat y_M$.
We further define
$f^\prime:=f\times *$, $f_U^\prime:=f_U\times*$, $f_V^\prime:=f_V\times *$, and we let
$f^\prime_M$ be the concatenation of the homotopy
$f_M\times *$ with the loop $*\times l_M\circ p$, where $p\colon I\times M\to \Sigma^u M$ is the obvious projection.
Then we have
$$\int_{I\times M/M} f^{\prime,*}_M R(\hat y^\prime)=\int_{I\times M/M} f_M^* R(\hat y)+
\int_{I\times M/M}   l_M^* R(\hat y_M)\ .$$
The closed form
 $\int l_M^* R(\hat y_M)$ represents the class $-\ch(I(\hat m))$.
Therefore, if we replace $f_M$ by  $f_M^\prime$, then  we can improve (\ref{uoiqdqwdqwwqd})
to 
$$
\kappa-\int_{I\times M/M} f_M^* R(\hat y)=d\sigma
$$
for some form $\sigma\in \Omega^{k-2}(M,\tV)$.
If we change the forms $\rho_U$ and $\rho_V$ in the construction of $\kappa$ to
$ \rho_U-\sigma_{|U}$ and $ \rho_V-\sigma_{|V}$, then with the resulting new choice for $\kappa$ we get
$$
\kappa-\int_{I\times M/M} f_M^* R(\hat y)=0\ .
$$
We form the loops
$l_V:= f_{M|V}^{op}\sharp f_V  \colon S^1\times V\to N$ and $l_U:=f_{M|U}^{op}\sharp f_U\colon S^1\times U\to N$.
We define
$$\hat u:=\int l_U^*\hat y +a(\rho_U)\ ,\quad \hat v:=\int l_V^* \hat y+a(\rho_V)\ .$$
Then we get
\begin{eqnarray*}
R(\hat u)&=&\int  l_U^* R(\hat y)+ d\rho_U\\&=&
 \int_{I\times U/U} f_U^* R(\hat y)-\int_{I\times U/U} f_{M|U}^* R(\hat y)+d\rho_U\\&=&\kappa_{|U}-\int_{I\times U/U} f_{M|U}^* R(\hat y)\\&=&0\\
R(\hat v)&=&0\ .
\end{eqnarray*}
Hence
$$\hat u\in \hat E^{k-1}_{flat}(U)\ ,\quad \hat v\in \hat E^{k-1}_{flat}(V)\ .$$
We now verify that
$$\hat v_{|A}-\hat u_{|A}=\hat x\ .$$
We have
$$\hat v_{|A}-\hat u_{|A}=\left(\int l_V^*\hat y\right)_{|A}- \left(\int l_U^* \hat y\right)_{|A}+a(\rho)=\left(\int l_V^*\hat y\right)_{|A}- \left(\int l_U^* \hat y\right)_{|A}+ \hat x-\int l^*\hat y\ .$$
It therefore suffices to show the following Lemma.
\begin{lem}
\begin{equation}\label{udiwqdqwdqwd}
\left(\int l_V^*\hat y\right)_{|A}- \left(\int l_U^* \hat y\right)_{|A}  =\int l^*\hat y\ .
\end{equation}
\end{lem}
\proof

 We choose an embedding of $S^1\tilde \vee S^1:=S^1\cup_{\{0\}=*}[0,1]\cup_{\{1\}=*}S^1$ into
$\R^2$ which is smooth on the two copies of $S^1$ and the interval and such that the interval intersects the circles \textcolor{black}{transversally}. 
A smooth function on $S^1\tilde \vee S^1$ is one which extends to a smooth function
on $\R^2$. We thus   have the notion of a  smooth map from $S^1\tilde \vee S^1$ to a manifold.
Moreover, a map   $W\to S^1\tilde \vee S^1$ from a manifold is smooth if the composition $W\to S^1\tilde \vee S^1\to \R^2$ is smooth.

We furthermore choose an open neighbourhood $W\subset \R^2$ of  $S^1\tilde \vee S^1$ which admits a smooth  projection $\pi\colon W\to 
S^1\tilde \vee S^1$ which is a homotopy equivalence.
We define $S^1\times A\tilde \vee_A S^1\times A:=(S^1\tilde \vee S^1\times A)$ and set 
$$\widetilde{l_U^{op}\vee l_V}\colon W\times A\stackrel{\pi\times \id}{\to} S^1\times A\tilde \vee_A S^1\times A\stackrel{l_U^{op}\tilde \vee l_V}{\to} N\ ,$$
where the smooth map $l_U^{op}\tilde \vee l_V$ maps the part $ [0,1]\times A$ to the base point and
is given by $l_U^{op}$ and $l_V$ on the left and right copies of $S^1\times A$, respectively.
We have a diagram
$$\xymatrix{&S^1\times A\sqcup S^1\times A\ar[d]^{j\sqcup k}&\\S^1\times A\ar@/^4cm/[rr]^e\ar[r]^s\ar@/^0.2cm/[rd]^{\tilde a}\ar@/_0.2cm/[dr]_{\tilde  b}&W\times A\ar[r]^{\widetilde{l_U^{op}\vee l_V}}\ar@/^0.5cm/[d]^a\ar@/_0.5cm/[d]^b&N\\&S^1\times A\ar@/^0.2cm/[ur]^{l_U^{op}}\ar@/_0.2cm/[ru]_{l_V}&}
$$
where $$a,b\colon W\times A\stackrel{\pi}{\to}S^1\times A\tilde \vee_A S^1\times A\to S^1\times A$$ are the  projections which contract the  left or right summand, respectively, and
$$j,k\colon S^1\times A\to S^1\times A\tilde \vee_A S^1\times A\to W\times A$$ are the embeddings of the left and right summand. The usual coproduct map $S^1\to S^1\vee S^1$ gives rise to a smooth  map
$$s\colon S^1\times A\to S^1\times A\tilde \vee_A S^1\times A\to W\times A\ .$$

We first observe that $j^*\oplus k^*\colon  E^*(W\times A)\to  E^*(S^1\times A)\oplus  E^*(S^1\times A)$ and $j^*\oplus k^*\colon \im(\pi^*)\to \Omega^*(S^1\times A,\tV)\oplus  \Omega^*(S^1\times A,\tV)$ are injective, where  
\begin{equation}\label{udqiwdqwdwqdqwdqwdwqd}
\im(\pi^*):=a^*\Omega^*(S^1\times A,\tV)+b^*\Omega^*(S^1\times A,\tV)\ .
\end{equation} Note that \textcolor{black}{the} definition of $\im(\pi^*)$ is a
slight abuse of notation. 
Next we show that $j^*\oplus k^*\colon \hat E^*(\textcolor{black}{W}
\times A)\to \hat E^*(S^1\times A)\oplus \hat E^*(S^1\times A)$ possesses a certain injectivity, too. 
Let $\hat r\in \hat E^*(W)$ be such that $R(\hat r)\in \im(\pi^*)$.
If $(j^*\oplus k^*)(\hat r)=0$, then $R(\hat r)=0$ and $I(\hat r)=0$. Therefore we can assume that
$\hat r=a(\rho)$ for some $\rho\in \Omega^{*-1}_{cl}(W)$.
Since $j^*a(\rho)=0$ and $k^*a(\rho)=0$ there exist classes
$s,t\in E^{*-1}(S^1\times A)$ such that
$\ch(s)=\Rham(j^*\rho)$ and $\ch(t)=\Rham(k^*\rho)$.
 Let $i\colon A\to S^1\times A$ be induced by the base point of $S^1$. 
Since $j\circ i$ is homotopic to $k\circ i$ we have  $i^*j^*\Rham(\rho)=i^*k^*\Rham(\rho)$. We therefore can in addition assume after modifying e.g. $t$ by a torsion class coming from $A$ that
$i^*s=i^*t$. But then there exists a class
$w\in E^{*-1}(W\times A)$ such that $j^*w=s$ and $k^*w=t$.
It follows that
$j^*\ch(w)=j^*\Rham(\rho)$ and $k^*\ch(w)=k^*\Rham(\rho)$.
This implies that
$\ch(w)-\Rham(\rho)=0$ and hence $a(\rho)=0$.

The composition $e$ is homotopic to the loop $l$ by a homotopy $H$. 
Indeed, the loop $e$ is the concatenation
$$f_{U|A}^{op}\sharp f_{M|A}^{op}\sharp f_{M|A}\sharp f_{|V|A}\ ,$$
where $F^{op}$ is the homotopy $F$ run in the opposite direction.
The homotopy $H$ to $f_{U|A}^{op}\sharp  f_{|V|A}$ can thus be arranged  symmetrically so that
$\int_{I\times S^1\times A/S^1\times A}H^*\omega=0$ for every $\omega\in \Omega(N)$.

Furthermore, the compositions $l_U^{op}\circ \tilde a,l_V\circ \tilde b$ are homotopic to
$l_U^{op}$ and $l_V$ by  homotopies of the form $G_U=g_U\times \id_A$ and $G_V=g_V\times \id_A$, where $g_U,g_V\colon I\times S^1\to S^1$. In the following we use the symbol $w$ in order to denote various constant maps to the base point of $N$.  
We have
\begin{equation}\label{uidqwdwqdqwdwqdfffef}
(\widetilde{l_U^{op}\vee l_V})^*\hat y-(l_U^{op}\circ a)^*\hat y-(l_V\circ b)^*\hat y=-w^*\hat y\ .
\end{equation}
Indeed, if we apply  $j^*$ to the left-hand side we get
$$j^*(\widetilde{l_U^{op}\vee l_V})^*\hat y-j^*(l_U^{op}\circ a)^*\hat y-j^*(l_V\circ b)^*\hat y=-
w^*\hat y\ .$$
Here we use
$$l_U^{op}\circ a=j^* (\widetilde{l_U^{op}\vee l_V}) \ ,\quad w=l_V\circ b\circ j\ .$$
Similarly we get
$$k^*(\widetilde{l_U^{op}\vee l_V})^*\hat y-k^*(l_U^{op}\circ a)^*\hat y-k^*(l_V\circ b)^*\hat y=-w^*\hat y$$
We now use the injectivity of $j^*\oplus k^*$. Note that
the curvature of both sides of (\ref{uidqwdwqdqwdwqdfffef}) are in $\im(\pi^*)$ as defined in (\ref{udqiwdqwdwqdqwdqwdwqd}).

Since the constant map $w\colon S^1\times A\to N$ factors over the projection $S^1\times A\to A$ we have $\int w^* \hat y=0$.
We calculate
$$e^*\hat y=s^*(\widetilde{l_U^{op}\vee l_V} )^*\hat y=l^*\hat y + a(\int_{I\times S^1\times A/S^1\times A} H^*R(\hat y))=l^*\hat y  \ .$$
Furthermore,
$$e^*\hat y=s^*(l_U^{op}\circ a)^*\hat y+s^*(l_V\circ b)^*\hat y-w^*\hat y=\tilde b^*l_V^*\hat y -\tilde a^*l_U^*\hat y-w^*\hat y\ ,$$
hence
$$e^*\hat y=l_V^*\hat y+l_U^{op,*}\hat y+a(\int_{I\times S^1\times A/S^1\times A}  G_V^*R(\hat y))
-a(\int_{I\times S^1\times A/S^1\times A}  G_U^*R(\hat y))-w^*\hat y\ .$$
We now apply $\int$ and observe using $\int l_U^{op,*}\hat y =-\int l_U^*\hat y$  and $\int w^*\hat y=0$ that
it suffices to show that
$$\int \int_{I\times S^1\times A/S^1\times A}  G_V^*R(\hat y)=0\ ,\quad 
\int\int_{I\times S^1\times A/S^1\times A}  G_U^*R(\hat y) =0\ .$$
Because of the special form of the homotopies $G_U$ and $G_V$ these integrals indeed vanish. This finishes the verification of (\ref{udiwqdqwdqwd}). \hB

\item 
Let $\hat x\in \hat E^k_{flat}(M)$ be such that $c(\hat x)=0$.
This means that $\hat x_{|U}=0$ and $\hat x_{|V}=0$. Let $x:=I(\hat x)$.
We choose a based  manifold $N$, a class $y\in E^k(N)$ with trivial restriction to the base point,  and a smooth map
$f\colon M\to N$ such that $f^*y=x$ and there are zero homotopies
$f_U\colon I\times U\to N$, $f_V\colon I\times V\to N$. We further choose a smooth lift $\hat y$ with trivial restriction to the base point
and a form $\lambda\in \Omega^{k-1}(M,\tV)$ such that $f^*\hat y=\hat x+a(\lambda)$.
>From the homotopy formula Lemma \ref{udqwdqwdqw} we get
\begin{eqnarray*}
0&=&\hat x_{|U}=-a(\int_{ I\times U/U} f_U^*R(\hat y)) +a(\lambda_{|U})\\0&=&\hat x_{|V}=-a(\int_{I\times V/V} f_V^*R(\hat y)) +a(\lambda_{|V})\ .
\end{eqnarray*}
Hence there exists classes $\hat u\in \hat E^{k-1}(U)$ and $\hat v\in\hat  E^{k-1}(V)$ such that
\begin{equation}\label{zdqiwdqwd}
\int_{I\times V/V} f_V^*R(\hat y)+\lambda_{|V}=R(\hat v)\ ,\quad \int_{I\times U/U} f_U^*R(\hat y)+\lambda_{|U}=R(\hat u)\ .
\end{equation}
We now show that by modifying the choices of $f\colon M\to N$, $y$ and the homotopies
we can assume that $R(\hat u)$ and $R(\hat v)$ are exact. By Lemma \ref{duiqdwqdqwd}
we choose a based manifold $N_U$, a class
$y_U$ which vanishes on the base point, and a map
$l_U\colon \Sigma^u U\to N_U$ such that
$\sigma(I(\hat u))=-l_U^*y_U$. We can assume that $l_U$ is constant near the singular points
of the unreduced suspension. We adopt a similar choice for $V$. We further choose
smooth lifts $\hat y_U$ and $\hat y_V$ again vanishing on the base points. 
We consider 
$N^\prime:=N\times N_U\times N_V$ with the class
$\pr_N^*\hat y+\pr_{N_U}^*\hat y_U+\pr_{N_V}^*\hat y_V$ and the map
$f^\prime:=f\times *\times *$.
We further concatenate the homotopy 
$(f_U\times *\times *)$ with the loop  $*\times l_U\times *$ in order to get
a new zero homotopy   $f^\prime_U$ of $f^\prime_{|U}$. We define $f^\prime_V$ in a similar manner. We now observe that
$$\int_{I\times V/V} f^{\prime,*}_V R(\hat y^\prime)=\int_{I\times V/V} f_V^*R(\hat y)+\int_{I\times V/V} l_V^*R(\hat y_V)$$ and the closed form
$\int_{I\times V/V} l_V^*R(\hat y_V)$ is in the cohomology class of $-R(\hat v)$.
A similar calculation holds for $U$. 

If we replace the old choices by the new choices we can now improve
(\ref{zdqiwdqwd}) to 
\begin{equation}\label{zdqiwdqwd1}
\int_{I\times V/V} f_V^*R(\hat y)+\lambda_{|V}= -d\rho_V\ ,\quad \int_{I\times U/U} f_U^*R(\hat y)+\lambda_{|U}=-d\rho_U\ .
\end{equation}
We define $\rho:=\rho_{V|A}-\rho_{U|A}$ and set
$$\hat z:=\int l^*\hat y +a(\rho)\ .$$
We calculate
\begin{eqnarray*}
R(\hat z)&=&\int l^* R(\hat y)+d\rho\\
&=&\left(\int_{I\times V/V} f_V^*R(\hat y)\right)_{|A}+d\rho_{V|A}-\left(\int_{I\times U/U} f_V^*R(\hat  y)\right)_{|A}-d\rho_{U|A}\\&=&
-\lambda_{|V|A}+\lambda_{|U|A}\\&=&0 \ .
\end{eqnarray*}
Hence $\hat z\in \hat E^{k-1}_{flat}(A)$. Furthermore, if we construct $\hat \partial \hat z$
using the choices fixed above we get
$$\kappa_{|U}=\int_{I\times U/U} f_U^*R(\hat  y)+d\rho_U=-\lambda_{|U}\ ,\quad  \kappa_{|V}=-\lambda_V\ .$$
This gives
$$\hat \partial \hat z=f^*\hat y-a(\lambda)=\hat x\ .$$
\item
Finally we show exactness at $\hat E^k_{flat}(U)\oplus \hat E^k_{flat}(V)$.
Let $$\hat u\in \hat E^k_{flat}(U)\ ,\quad \hat v\in  \hat E^k_{flat}(V)$$ 
be such that $\hat u_{|A}=\hat v_{|A}$. 
Let $u:=I(\hat u)$ and $v:=I(\hat v)$.
By Lemma \ref{7eifdwqedqfqd}
we choose a manifold $N$ with a class $y\in E^k(N)$, smooth maps
 $f\colon U\to N$, $g\colon V\to N$ such that $f^*y=u$ and $g^*y=v$, and there is a homotopy
$f_{|A}\sim g_{|A}$ which we denote by $h$.

In a  first step  we  show that we can choose a map $e\colon M\to N$ such that
\begin{equation}\label{uasidasd}
e^*_{|U}y=u\ ,\quad e^*_{|V}y=v\ .
\end{equation} In fact we can define
\begin{equation}\label{uasidasd1}e(m):=\left(\begin{array}{cc}f(m)&m\in M\setminus V\\
h(\frac{\chi(m)+1}{2},m)&m\in A\\
g(m)&m\in M\setminus U
\end{array}\right)\ .\end{equation} 
The relations (\ref{uasidasd}) hold true since there are homotopies $e_{|U}\sim f$ and $e_{|V}\sim g$.

We choose $\alpha_U\in \Omega^{k-1}(U,\tV)$ and   $\alpha_V\in \Omega^{k-1}(V,\tV)$ such that 
$$e_{|U}^*\hat y+a(\alpha_U)=\hat u\ ,\quad e_{|V}^*\hat y+a(\alpha_V)=\hat v\ .$$
We then have
$a(\alpha_{V|A}-\alpha_{U|A})=0$ so that
$\alpha_{V|A}-\alpha_{U|A}=R(\hat w)$ for some $\hat w\in \hat E^{k-1}(A)$.
In the following we show that by modifying the homotopy $h$ we can 
assume that $R(\hat w)$ is exact.

 Using Lemma \ref{duiqdwqdqwd} we choose a based manifold
$N^\prime$, a class $y^\prime\in E^{k-1}(N^\prime)$, and a map $l\colon \Sigma^u A\to N^\prime$ such that $l^*y^\prime=-\sigma (I(\hat w))$.
We can assume that $l$ maps $(t,a)\in \Sigma^u A$ to the base point  if
$t\in [0,1/4]$ or $t\in [3/4,1]$.

Without loss of generality we can assume that the homotopy $h\colon I\times A\to N$ is
constant on the part $[1/4,3/4]\times A$.
We now replace
$N$ by $\tilde N:=N\times N^\prime$, $y$ by $\tilde y=\pr_N^*y+\pr_{N^\prime}^*y^\prime$, $f$ by $\tilde f:=f\times *$, $g$ by $\tilde g:=g\times *$ and
$h$ by $\tilde h\colon I\times A\to \tilde N$ given by
$$\tilde h(t,a):=\left(\begin{array}{cc}
(h(t,a),*)&t\in [0,1/4]\\
(h(1/2,a),l(t,a))&t\in [1/4,3/4]\\
(h(t,a),*)&t\in [3/4,1]
\end{array}\right)
\ .$$
Let $\tilde e\colon M\to \tilde N$ be the resulting map (\ref{uasidasd1}).
Note that there is a homotopy $d_U$ from
$\tilde e_{|U}$ to $e_{|U}\times *$ and a similar homotopy $d_V$ from $\tilde e_{|V}$ to $e_{|V}\times *$.
We furthermore set $\hat {\tilde y}:=\pr_N^*\hat y+\pr_{N^\prime}^*\hat y^\prime$
for some smooth lift $\hat y^\prime\in \hat E^{k-1}(N^\prime)$ which we arrange such that $\hat y^\prime_{|*}=0$.
Since $(e_{|U}\times *)^*\hat{\tilde y}=e_{|U}^*\hat y$ we can choose
$$\tilde \alpha_U:=\alpha_U+\int_{I\times U/U} d_U^*R(\hat {\tilde y})\ , \quad \tilde \alpha_V:=\alpha_V+\int_{I\times V/V} d_V^*R(\hat {\tilde y}) .$$
We get
$$\tilde \alpha_{V|A}-\tilde \alpha_{U|A}=R(\hat w)+\int  L^*R(\hat {\tilde y})\ ,$$
where
$L\colon S^1\times A\to \tilde N$ is the loop $L:=d_U^{op}\sharp d_V$.
Note that we can choose $d_U$ and $d_V$ such that
$\pr_N\circ L$ is constant. Therefore
$$\int  L^*R(\hat {\tilde y})=\int  L^*\pr_{N^\prime}^* R( \hat y^\prime)\ .$$
We get for the cohomology classes
$$\Rham(\int  L^*\pr_{N^\prime}^* R(\hat  y^\prime))=\int \ch(\sigma(I(\hat w)))=-\ch(\hat w)\ .$$
If we replace the objects without a tilde by the objects with the tilde decoration, then
we can assume that $\alpha_{V|A}-\alpha_{U|A}=d\sigma$.
We choose $\sigma_U\in \Omega^{k-2}(U,\tV)$ and $\sigma_V\in \Omega^{k-2}(V,\tV)$ such that
$\sigma_{V|A}-\sigma_{U|A}=\sigma$.
Then we replace $\alpha_U$ by $\alpha_U-d\sigma_U$ and $\alpha_V$ by $\alpha_V-d\sigma_V$. After these changes we can assume that
$\alpha_{V|A}=\alpha_{U|A}$, hence $\alpha_U$ and $\alpha_V$ are restrictions of a global $\alpha\in \Omega^{k-1}(M,\tV)$.
We define $$\hat x:=e^*\hat y+ a(\alpha)\ .$$
Then $$\hat u=\hat x_{|U}\ ,\quad \hat v=\hat x_{|V}\ .$$
It also follows that $\hat x\in E^k_{flat}(M)$. This provides the required preimage of the sum
$\hat u\oplus \hat v$.
\end{enumerate}
\hB

We now have a homotopy invariant functor $\hat E^*_{flat}$ defined on smooth manifolds with a natural Mayer-Vietoris sequence. It gives rise to a similar functor on the category of finite-dimensional countable $CW$-complexes by the following proposition. 
 \begin{prop}\label{uidwqdqwd}
A  homotopy invariant functor $$H\colon \{\mbox{\tt smooth manifolds}\}\to \{\Z-\mbox{\tt graded abelian groups}\}$$
  with a natural Mayer-Vietoris sequence gives extends uniquely to a homotopy invariant functor
$$h\colon \{\mbox{\tt finite-dimensional countable $CW$-complexes}\}\to  \{\Z-\mbox{\tt graded abelian groups}\}$$  with a natural Mayer-Vietoris sequence.
\end{prop}
\proof
For a proof we refer to 
\cite{ks}. It uses the fact that diagrams of maps between countable
finite-dimensional $CW$-complexes can be approximated up to homotopy by
corresponding diagrams of manifolds. \hB 
Also the following lemma is well-known.
\begin{lem}\label{uoqwdqwdwqd}
A functor $$h\colon \{\mbox{\tt finite $CW$-complexes}\}\to   \{\Z-\mbox{\tt graded abelian groups}\}$$  with a natural Mayer-Vietoris sequence gives rise to a reduced cohomology theory  $\tilde h$ on the category of pointed finite $CW$-complexes.
\end{lem}
\proof
For a finite pointed $CW$-complex $X$ we define 
$$\tilde h^*(X):=\ker(h^*(X)\to h^*(*))\ .$$
To each map $f\colon X\to Y$ of pointed $CW$-complexes we get an induced map
$f^*\colon \tilde h^*(Y)\to \tilde h^*(X)$.
Let $C(X):=[0,1]\times X/\{1\}\times X$ be the cone over $X$ with its natural $CW$-structure. Then
we can write the unreduced suspension as
\begin{equation}\label{udioqwdwqd}
\Sigma^uX=C(X)\cup_X C(X)\ .
\end{equation} 
The suspension isomorphism 
$$\sigma\colon \tilde h^*(X)\to \tilde h^*(\Sigma^u(X))$$
is given by the boundary operator in the Mayer-Vietoris sequence
associated to the decomposed $CW$-complex (\ref{udioqwdwqd}).
It is obviously natural.
Finally, for each subcomplex $A\subseteq X$ the mapping cone sequence
$$A\to X\to X\cup_AC(A)$$ gives rise to an exact sequence
$$\tilde h^*( X\cup_AC(A))\to \tilde h^*(X)\to \tilde h^*(A)\ .$$
Indeed, this is a part of the Mayer-Vietoris sequence for the decomposition
$X\cup_AC(A)$ since $\tilde h^*(CA)=0$.
\hB

If $\tilde h$ is a reduced cohomology theory on the category of pointed finite $CW$-complexes, then by \cite[Thm. 9.27]{MR1886843} there exists a spectrum
$\bh$ which represents $\tilde h$. The isomorphism class of this spectrum is well-defined.
Furthermore by \cite[Thm. 9.27]{MR1886843}, a natural transformation $\tilde h\to \tilde h^\prime$ of reduced cohomology theories on finite $CW$-complexes can be represented by a map of spectra $\bh\to \bh^\prime$, which might be not uniquely determined.

If we apply these topological results Proposition \ref{uidwqdqwd}, Lemma \ref{uoqwdqwdwqd}
to the flat theory $\hat E_{flat}^*$, then we get a reduced cohomology theory
$\tilde U^{*+1}$ on the category of pointed finite $CW$-complexes which we can represent by a spectrum $\bU$
whose isomorphism class is well-defined.
Since every compact manifold has the structure of  a finite $CW$-complex we can restrict the
theory $\tilde U^*$ again to compact manifolds. We thus have shown:

\begin{theorem} \label{udidowqdiwqdwqd}
If $(\hat E,R,I,a,\int)$ is a smooth extension of $E$ with integration, then  $\hat E^*_{flat}$
has a natural long exact Mayer-Vietoris sequence. Its restriction to compact manifolds
is equivalent to the restriction  to compact manifolds of a generalized cohomology theory represented by a spectrum.
\end{theorem}

We now compare $\hat E_{flat}^{*-1}$ with $E\R/\Z^*$.
The natural transformation
$H^{*-1}(M;\tV)\to \hat E^{*}_{flat}(M)$ induced by $a$ gives a natural transformation
$\tilde E\R^{*}\to \tilde U^*$ which  can be represented by a map of spectra  $\bE\R\to \bU$.
 
We now consider the diagram of distinguished triangles in the stable homotopy category
$$\xymatrix{\bF ibre\ar[r]&\bE\R\ar[r]\ar@{=}[d]& \bU\ar[r]&\Sigma\bF ibre\\
\bE\ar@{-->}[u]\ar[r]\ar@{.>}[rru]&\bE\R\ar[r]&\bE\R/\Z\ar@{-->}[u]\ar[r]&\Sigma \bE\ar@{-->}[u]}\ .$$
The fact that
$$E^*(M)\stackrel{\ch}{\to} H^*(M;\tV)\stackrel{a}{\to} \hat E(M)$$ vanishes
implies that the dotted arrow is trivial. This gives the dashed factorization $\bE\R/\Z\to \bU$
which we extend to a map of triangles. Note that the dashed maps are not necessarily unique.

\begin{theorem}\label{udqiwdqwdqwddwqdqwd1}
Assume that $(\hat E,R,I,a \int)$ is a smooth extension of a generalized cohomology with integration.
If $E^*$ is torsion-free, then
there exists a natural isomorphism (not necessarily unique) of functors on compact manifolds
$E\R/\Z^*(M)\to \hat E^{*-1}_{flat}(M)$
so that
$$\xymatrix{H^{*-1}(M;\tV)\ar[r]^a\ar[d]^\cong&\hat E^{*}_{flat}(M)\ar[r]&E^*(M)\ar[d]\\
E\R^{*-1}(M)\ar[r]&E\R/\Z^{*-1}(M)\ar[u]^\cong\ar[r]&E^*(M)}$$
commutes
\end{theorem}
\proof
We know by Theorem \ref{udidowqdiwqdwqd} that there is a natural isomorphism
$ \hat E^{*}_{flat}(M)\cong U^{*-1}(M)$.
It suffices to check that the transformation $\bE\R/\Z\to \bU$ is an equivalence by working on the level of homotopy groups. In other words, we must show that it induces an isomorphism on coefficients. We know that
$$\coker(a\colon E\R^k(M)\to \hat E^{k+1}_{flat}(M))\cong E^{k+1}_{tors}(M),\quad
\textcolor{black}{\ker(a\colon E\R^n(M)\to \hat E^{n+1}_{flat}(M))=\im(\ch)}\ .$$
\textcolor{black}{Therefore} we have a morphism of exact sequences
$$\xymatrix{E^k\ar[r]^\ch\ar@{=}[d]&E\R^k\ar[r]\ar@{=}[d]& U^k\ar[r]&E^{k+1}_{tors}\ar[r]&0\\
E^k\ar[r]^\ch&E\R^k\ar[r]&E\R/\Z^k\ar[r]\ar[u]&E^{k+1}_{tors}\ar@{=}[u]\ar[r]&0}\ .$$
If   $E^{k+1}_{tors}=0$, then the morphism $E\R/\Z^k\to U^k$ is an isomorphism
\textcolor{black}{by the Five-lemma}.
\hB

\section{Absence of Phantoms}

In the proof of Proposition \ref{udqiwdwqdqwd} we have used the fact that
certain generalized cohomology groups are \textcolor{black}{free of
  phantoms}. The absence of phantoms might be interesting also in other cases
where an approximation of an infinite loop space by manifolds is
invoked. Therefore we add this section. The results are probably well known,
but we couldn't find appropriate references.

In the following we assume that $E$ is a cohomology theory represented by a
commutative ring spectrum $\bE$. Let $Z$ be a $CW$-complex. 
\begin{ddd}\label{dqwudqwdqw}
We define the subspace of  \textbf{phantom classes} $E^*_{phantom}(Z)\subseteq E^*(Z)$ to be the subspace of all classes  $\phi\in E^*(Z)$ such that $f^*\phi=0$ for all maps $f\colon X\to Z$ and finite complexes $X$. 
\end{ddd}
In the following we discuss various conditions implying the absence of non-trivial phantom classes.


\begin{prop}\label{pr0}
If $E_*(Z)$ is a free $E^*$-module, then $E^*_{phantom}(Z)\cong 0$.
\end{prop}
\proof
We equip $E^*(Z)$ with the  profinite filtration topology
induced by the submodules $F^aE^*(Z):=\ker\left(E^*(Z)\to E^*(Z_a)\right)$, where
$(Z_a)$ is the system of all finite subcomplexes of $Z$.
On the other hand we equip
$DE_*(Z):=\Hom_{E^*}(E_*(Z),E^*)$ with the dual finite topology
generated by the submodules
$\ker(D E_*(Z)\to DL_b)$, where $(L_b)$ runs over the system of all finitely generated submodules of $E_*(Z)$. With this topology the $E^*$-module $DE_*(Z)$ is complete and Hausdorff.
By \cite[Thm. 4.14]{MR1361899} the evaluation
$E^*(Z)\otimes E_*(Z)\to E^*$ induces a topological isomorphism
$E^*(Z)\to D E_*(Z)$. The fact that the profinite filtration topology on $E^*(Z)$
is Hausdorff is equivalent to the absence of phantom classes. \hB 

\begin{lem}\label{zueweq}
$E^*_{Phantom}(\bE_k)\cong 0$ for $\bE=\bMU$ or even $k$ and $\bE=\bK$.
\end{lem}
\proof
The cohomology theories $MU^*$ and $K^*$ are represented by ring spectra.
We first consider the case $MU$.
In \cite[Sec. 4]{MR0356030} it \textcolor{black}{is} shown that $MU_*(\bMU_k)$ is a free
$MU^*$-module (it actually has been calculated completely).
We can therefore apply Proposition  \ref{pr0}.

For $K$-theory we we first note that
$\bK_0\cong \Z\times BU$, and that
$K_*(BU,\Z)$ is a free $\Z$-module
on even generators (\cite[Prop. 4.3.3 (d)]{MR1407034}).  
This implies that $K_*(\bK_0)$ is a free $K^*$-module, and we can again apply
Proposition \ref{pr0}. \hB

\begin{prop}\label{pr1}
If $E_*(Z)$ is a free $E^*$-module, then $E^*_{Phantom}(Z\times Z)\cong 0$.
\end{prop}
\proof
If $E_*(Z)$ is a free $E^*$-module, then so is $E_*(Z\times Z)$. In fact,
for every complex $X$ we have the K\"unneth isomorphism
$$E_*(Z)\otimes_{E^*}E_*(X)\stackrel{\sim}{\to} E_*(Z\times X) \ .$$
This follows from the usual observation that 
$E_*(Z)\otimes_{E^*}E_*(\dots)  \to E_*(Z\times \dots)$ is a natural transformation of homology theories
which coincide on the point. 
Finally, we use the fact that the tensor product of two free modules is again
free.\hB

\begin{kor}\label{hdwqdqwdwq}
If $k$ is even, then
$$E^{*}_{Phantom}(\bE_k\times \bE_k)\cong 0$$
holds true for
$\bE\in \{\bMU,\bK\}$.
\end{kor}
%

If $X\mapsto E^*(X)$ is a cohomology theory represented by a spectrum $\bE$, then let $X\mapsto E_*(X)$ denote the associated homology theory.
 We define
$$E_\R^*(X):=\Hom_{\Ab}(E_*(X),\R)\ ,\quad E_{\R/\Z}^*(X):=\Hom_{\Ab}(E_*(X),\R/\Z)\ .$$ 
Since $\R$ and $\R/\Z$ are injective abelian groups these constructions define
cohomology theories on the category of all topological
  spaces. Since they satisfy in addition the wedge axiom 
they can be represented by spectra which we denote by
$\bE_\R$ and $\bE_{\R/\Z}$.


\begin{lem}\label{dqduqwdwqdwq}
For every $CW$-complex $X$ we have  $E^*_{\R,phantom}(X)\cong 0$ and
$E^*_{\R/\Z,phantom}(X)\cong 0$.
\end{lem}
\proof
 Let us discuss the case of $E_\R$.  The case of $E_{\R/\Z}$ is similar.
It suffices to show that
$$E^*_{\R}(X)\cong \lim_a E^*_{\R}(X_a)\ ,$$
where $(X_a)$ is the system of finite subcomplexes of $X$.
We have $X\cong \colim_a\:  X_a$.
Since homology is compatible with colimits and
the $\Hom_\Ab(\dots,\R)$-functor turns colimits in the first argument into limits we get (compare
\cite[A.9]{frlo})
\begin{eqnarray*}
E^*_{\R}(X)&=&\Hom_\Ab(E_*(X),\R)\\
&\cong&\Hom_\Ab(E_*(\colim_a \: X_a),\R)\\
&\cong&\Hom_{\Ab}(\colim_a\:  E_*(X_a),\R)\\
&\cong&\lim_a\:  \Hom_{\Ab}(E_*(X_a),\R)\ . 
\end{eqnarray*}

\hB

%

The projection $\R\to \R/\Z$ induces a natural transformation of
cohomology theories $E_\R^*(X)\to E^*_{\R/\Z}(X)$. It is given by a morphism of representing spectra $\bE_\R\to \bE_{\R/\Z}$.
\begin{ddd}
The Andersen dual $D(E)$ of the cohomology theory $E$ is defined as the cohomology
theory represented by the spectrum $D(\bE)$ obtained by the extension of the map
$\bE_\R\to \bE_{\R/\Z}$ to a distinguished triangle in the stable homotopy catgeory
$$D(\bE)\to \bE_\R\to \bE_{\R/\Z}\to \Sigma D(\bE)\ .$$
\end{ddd}
 In \cite[p.~244]{MR704613} a morphism of distinguished triangles in the stable homotopy category 
$$\xymatrix{D(\bE)\ar[r]\ar@{=}[d]&D(\bE)\R\ar[r]\ar[d]^\cong& D(\bE)\R/\Z\ar[d]\ar[d]^\cong\ar[r]&\Sigma D(\bE)\ar[d]_{-1}^\cong \\D(\bE)\ar[r]&\bE_\R\ar[r]&\bE_{\R/\Z}\ar[r]&\Sigma D(\bE)}
$$
has been constructed so that the vertical maps are  equivalences. 

We now assume that $E^k$ is finitely generated for every $k\in \Z$.
Since $D$ is a duality on cohomology theories with finitely generated
coefficients \cite[Thm.~33]{MR704613}
we get by inserting $D(\bE)$ in place of $\bE$ and using $D(D(\bE))\cong \bE$ that
 $$\xymatrix{\bE\ar[r]\ar@{=}[d]&\bE\R\ar[r]\ar[d]^\cong&
   \bE\R/\Z\ar[d]^\cong\ar[r]&\Sigma \bE\ar[d]^\cong_{-1}\\
\bE\ar[r]&D(\bE)_\R\ar[r]&D(\bE)_{\R/\Z}\ar[r]&\bE}\ , $$
 i.e.~\textcolor{black}{we get} in particular
isomorphisms 
\begin{equation}\label{zdqwduqwdd}
\bE\R\cong  D(\bE)_\R ,\quad \bE\R/\Z\cong D(\bE)_{\R/\Z}\ .
\end{equation}

Combining (\ref{zdqwduqwdd}) with Lemma \ref{dqduqwdwqdwq}
(applied to $D(E)$ in the place of $E$)
we get
\begin{kor}
If $E$ is a cohomology theory represented by a spectrum such that
$E^k$ is finitely generated for all $k\in \Z$, then for every $CW$-complex $X$
we have
$$E\R^*_{Phantom}(X)=0\ ,\quad E\R/\Z^*_{Phantom}(X)=0\ .$$
\end{kor}

\begin{kor}\label{uoewqeqwe}
We have $K\R/\Z^0_{Phantom}(\bK_1\times \bK_1)=0$.
\end{kor}

\end{document}